\documentclass[12pt,reqno]{amsart}

\usepackage{color}   
\usepackage{hyperref}
\hypersetup{
  colorlinks=true, 
   linkcolor=red,  
}
 
 \newtheorem{theorem}{Theorem}
 \newtheorem{lemma}{Lemma}

 \newcommand{\eq}{\begin{equation}}
 \newcommand{\en}{\end{equation}}
 
 \def\proof{\noindent{\bf Proof\ \ }}
 \def\indist{\stackrel{\mbox{\small d}}{=}}
 \def\inprob{\mbox{$\rightarrow_P\ $}}
 \def\todist{\Rightarrow}
\def\bothsides{\stackrel{\mbox{\scriptsize $\smile$}}{\mbox{\scriptsize 
 $\frown$}}}
 
 \def\qed{\mbox{\rule{0.5em}{0.5em}}}
 
 \def\B#1{{\Bbb #1}}
 
 \def\BP{{\Bbb P}}
 \def\BR{{\Bbb R}}
 \def\BN{{\Bbb N}}
 \def\BE{{\Bbb E}}
 \def\BZ{{\Bbb Z}}
 
 \def\b#1{\mbox{\boldmath $#1$}}
 \def\ba{\mbox{\boldmath $a$}}

 \def\bu{\mbox{\boldmath $u$}}
 \def\bC{\mbox{\boldmath $C$}}
\def\bD{\mbox{\boldmath $D$}}
 \def\bZ{\mbox{\boldmath $Z$}}
 \def\bone{\mbox{\bf 1}}
 
 \def\smb#1{\mbox{\boldmath $#1$}}

 \def\TVone{{\rm Arratia and Tavar\'e (1992a)}}
 \def\TVtwo{{\rm Arratia and Tavar\'e (1992b)}}
 \def\TVthree{{\rm Arratia, Barbour and Tavar\'e (1992)}}
 \def\d{\displaystyle}
 \def\no{\noindent}
 
\begin{document} 

 \title[IPARCS (Adv. Math. 1994 pages 90 -- 154)]{Independent Process Approximations for Random Combinatorial
 Structures}
 \author{Richard Arratia and Simon Tavar\'e
  } 
 \thanks{Department of Mathematics,
 University of Southern California, Los Angeles, CA 90089-1113. 
 The authors were supported in part by NSF grant
 DMS90-05833. We thank Andrew Barbour, B\'ela Bollob\'as,
 Persi Diaconis, Jennie Hansen, Lars Holst, Jim Pitman and Gian-Carlo Rota for 
 helpful comments on earlier drafts of this paper.}
 \date{
{\sl Advances in Mathematics}, {\bf 43}, 000-000, 1994.\\ 
\vskip 2pc Dedication: To the memory of Mark
Kac,\\ who sought
out independence in \\ analysis and number theory}
 \maketitle
 
 \newpage
 \begin{center}
 {\bf Abstract}
 \end{center}
 
 {\small 

 Many random combinatorial objects have a component structure whose joint
 distribution is  equal to that of a process of mutually
 independent random variables, conditioned on the value of a weighted
 sum of the variables.
  It is interesting to compare the combinatorial structure directly
 to the independent discrete process, without renormalizing. The quality of
 approximation can often be conveniently quantified in terms of total variation
 distance, for functionals which observe part, but not all, of the
 combinatorial and independent processes.  
 
Among the examples are combinatorial assemblies (e.g. permutations, random
mapping functions, and partitions of a set), multisets (e.g. polynomials over a
finite field, mapping patterns and partitions of an integer), and selections
(e.g. partitions of an integer into distinct parts, and square-free
polynomials over finite fields).
  
 We consider issues common to all the above examples, including equalities 
 and upper bounds for total variation distances, existence of limiting
processes, heuristics for good 
 approximations, the relation to standard generating functions,
 moment formulas and recursions for computing densities, refinement to the
 process which counts the number of parts of each possible type, the effect of 
 further conditioning on events of moderate probability,
  large deviation theory and nonuniform 
 measures on combinatorial objects, and the possibility of getting useful
 results by overpowering the conditioning.

 \vskip 3pc
 
 \no Running Head: Independent Process Approximations
 
 \no Keywords: Poisson process, total variation, partitions,
permutations, random mappings, random polynomials, assemblies,
multisets, selections, Ewens sampling formula,
  generating functions, conditional limit theorems, sufficient
 statistics, large deviation theory
 
 \no MSC classification: 60C05, 60F17, 05A05, 05A16
 }
 \newpage

\tableofcontents 
 \section{Introduction}\label{sect1}

We consider random combinatorial objects which can be described in
terms of their component structure.  For an object of weight $n$, denote
the component structure by 
$$
\b{C} \equiv \b{C} (n) \equiv (C_1(n),C_2(n),\ldots,C_n(n)),
$$
where $C_i \equiv C_i(n)$ is the number of components of size $i$. Since
$i C_i$ is the total weight in components of size $i$, we have 
$$
C_1 + 2C_2 + \cdots + n C_n =n.
$$
For each fixed $n$, by choosing an object of weight $n$ at random, with
all possibilities equally likely, we view $\b{C} (n)$
as a $\B{Z}_+^n$-valued stochastic process, whose coordinates $C_i(n),
i=1,\ldots,n$, are {\it dependent}, nonnegative integer--valued 
random variables.  
This paper considers combinatorial objects for
which the joint distribution of $\b{C}(n)$ can be expressed as the joint
distribution of {\it independent} random variables $Z_1, Z_2, \ldots, Z_n $
conditioned on the value of a particular weighted sum.

    There are at least three broad classes of combinatorial structures
which have this description in terms of conditioning an independent
process.  The first class is assemblies of labeled structures on
$[n] \equiv \{1,2,\ldots,n \}$; see Joyal (1981). This class includes  
permutations, decomposed into cycles; mappings, decomposed into
connected components; graphs, decomposed into connected components, 
and partitions of a finite set.  The second class
is multisets, i.e. unordered samples taken with replacement. This class
includes 
partitions of an integer; random mapping patterns; and monic polynomials
over a finite field, decomposed into monic irreducible factors.  The
third class is selections, i.e. unordered samples taken without replacement, 
including partitions of an integer into parts of distinct sizes, and
square-free polynomials.

    The detailed description of any of the above examples is given in
terms of a sequence of nonnegative integers $m_1, m_2, \ldots $ \  .
 For assemblies, let $m_i$ be the number of labelled structures 
on a set of size $i$, for $i=1,2,\ldots $ ;
 permutations have $m_i =(i-1)!$,  mappings have
$ m_i = (i-1)!(1 + i + i^2 /2 +... + i^{i-1}/(i-1)!)$, and partitions
of a set have  $m_i=1$.
For multisets and selections, let $m_i$ be the number of objects of
weight $i$; partitions of an integer have $m_i=1$,
and the factorizations of monic polynomials over a finite field have $m_i$
equal to the number of monic, irreducible polynomials of degree $i$.

For $\b{a} \equiv (a_1,a_2, \ldots, a_n) \in \B{Z}^n_+$, consider the
number  $N(n,\b{a})$ of
objects of total weight $n$, having $a_i$ components of size $i$, for $i
= 1$ to $n$.  For assemblies, the generalization  of Cauchy's formula for
permutations is the enumeration
\begin{eqnarray}
N(n, \b{a}) & \equiv &  | \{ \mbox{assemblies on $[n]$}: \b{C} = \b{a} \} |
\nonumber \\
& = &
 \bone(a_1+2 a_2 + \cdots + n a_n =n) 
\   n! \ \prod_1^n \frac{m_i^{a_i}}{(i!)^{a_i} \  a_i !}.
\label{Nassembly}
\end{eqnarray}
For multisets, 
\begin{eqnarray}
N(n, \b{a}) & \equiv &  | \{ \mbox{multisets of weight $n$}: \b{C} = \b{a}
\} |
 \nonumber \\
& = &
 \bone(a_1+2 a_2 + \cdots + n a_n =n) 
\    \ \prod_1^n {m_i +a_i -1 \choose a_i }.
\label{Nmultiset}
\end{eqnarray}
For selections, 
\begin{eqnarray}
N(n, \b{a}) & \equiv &  | \{ \mbox{selections of weight $n$}: \b{C} = \b{a}
\} |
  \nonumber \\
& = &
 \bone(a_1+2 a_2 + \cdots + n a_n =n) 
\    \ \prod_1^n {m_i  \choose a_i }.
\label{Nselection}
\end{eqnarray}

Let $p(n)$ denote the total number of structures of weight $n$, to wit 
\eq\label{def p(n)}
p(n) = \sum_{\b{a} \in \B{Z}^n_+} N(n,\b{a}).
\en
For permutations, $p(n)=n!$; for mappings, $p(n)=n^n$; for graphs, $p(n)
= 2^{n \choose 2}$; for partitions of
a set $p(n) = B_n$, the Bell number; for partitions of an integer,
$p(n)$  is the standard notation; and for monic polynomials over a field
with $q$ elements, $p(n)=q^n$.

A {\it random structure} is understood as follows. Fix a constant $n$, and
choose one of the $p(n)$ structures at random, with each possibility
equally  likely.  This makes $\b{C}(n)$ a
stochastic process with values in $\b{Z}_+^n$, whose distribution is
determined by 
\eq\label{combdist.old}
\B{P}(\b{C}(n)=\b{a}) \equiv \frac{N(n,\b{a})}{p(n)},\ \b{a} \in
\B{Z}_+^n.
\en

In Section \ref{sect2} below, we show that there are independent random
variables $Z_1, Z_2, \ldots$ such that
the combinatorial
distribution (\ref{combdist.old}) is equal to the joint distribution of
$(Z_1,Z_2,\ldots,Z_n)$ conditional on the event
$\{T_n =n \}$, where 
$$
T_n \equiv Z_1+2Z_2+\cdots+nZ_n.
$$
Explicitly, for all 
$\b{a} \in \B{Z}_+^n$
\eq\label{equaldist}
\B{P}(\b{C}(n)=\b{a}) = \B{P}\left( (Z_1,Z_2,\ldots,Z_n)=\b{a} \left|
 T_n=n
\right. \right).
\en

The distribution of the degrees of the factors of the characteristic 
polynomial of a uniformly chosen random
matrix  over a finite field can also be expressed as in
(\ref{equaldist}) (cf.
Hansen and Schmutz (1994)), but it is not an example of an assembly,
multiset, or selection.

  It is fruitful to compare the combinatorial structure directly
 to the independent discrete process, without renormalizing. The quality of
 approximation can  be usefully quantified in terms of total variation
 distance between the restrictions of the dependent and independent
processes to a subset of the possible coordinates. We carry this out in
Section \ref{sect3}. 
Bounds and limit theorems for natural functionals which depend on the
coordinates, albeit weakly on those outside a subset, are then easily
obtained as corollaries. For examples of this in the context of random
polynomials over finite fields, and random
permutations and random mappings, see Arratia, Barbour, and Tavar\'e
(1993), and Arratia and Tavar\'e (1992b).

The comparison of combinatorial structures to independent processes,
with and without further conditioning, has
a long history.  Perhaps the best known example is the representation of
the multinomial distribution with parameters $n$ and $p_1,\ldots,p_k$ as
the joint law of independent Poisson random variables with means
$\lambda p_1, \ldots, \lambda p_k$, conditional on their sum being equal
to $n$.

Holst (1979a) provides an approach to urn models that unifies 
multinomial, hypergeometric and P\'olya sampling. The joint laws of the
dependent counts of the different types sampled are represented,
respectively, 
as the joint  distribution of independent Poisson,
negative binomial, and binomial random variables, conditioned
on their sum. See also Holst (1979b, 1981). The quality
of such approximations  is assessed using metrics, including the total
variation distance, by Stam (1978) and Diaconis and Freedman (1980).  

The books by Kolchin, Sevast'yanov, and
Chistyakov (1978) and Kolchin (1986) use the representation of
combinatorial structures, including random permutations and random
mappings, in terms of independently
distributed random variables, conditioned on the value of their sum.
However, the Kolchin technique requires that the independent variables be {\em
identically} distributed.  The
number of components $C_i$ of size $i$ is the number of random variables
which take on the value $i$.

Shepp and Lloyd (1966) study random permutations using
 a conditional relation almost identical
to (\ref{equaldist}), with $\BE Z_i=x^i/i$ and $x=x(n)$, 
except that they condition on $n$ being the value
of an infinite sum $Z_1+2Z_2+\cdots$, which of course entails that
$Z_{n+1}=Z_{n+2}=\cdots=0$, and requires $x<1$. Variants on the
 Shepp and Lloyd technique are discussed by Diaconis and Pitman (1986),
and are effectively exploited to prove functional
central limit theorems for two combinatorial assemblies by Hansen (1989,
1990), and used as a convenient tool for moment calculations by Watterson
(1974a) and Hansen
(1993). A related technique, coupled with an observation of Levin (1981), is
used by Fristedt (1992, 1993) to study random partitions of a set and
random partitions of an integer.

\subsection{Notation}

There are several types of asymptotic relation used in this paper. For
sequences $\{a_n\}$ and $\{b_n\}$, we write
 $a_n \sim b_n$ for the asymptotic relation
 $a_n/b_n \to 1$ as $n \to \infty$. We write
 $a_n \bothsides b_n$ if there are constants $0 < c_0 \leq c_1 <
\infty$ such that $c_0 b_n
\leq a_n \leq c_1 b_n$ for all sufficiently large $n$. 
We write $a_n \approx b_n$ to denote that $\log a_n \sim \log b_n$. 
Finally, we say that $a_n \doteq b_n$ if $a_n$ and
$b_n$ are approximately equal in some heuristic sense
deliberately left vague.

For $r \in \BZ_+ \equiv \{0,1,2,\ldots \}$,
we denote the rising factorial $y_{(r)}$ by $y_{(0)} = 1, 
y_{(r)} = y(y+1) \cdots (y+r-1)$  and the falling factorial $y_{[r]}$ by 
$y_{[0]} = 1, y_{[r]} = y(y-1) \cdots (y-r+1)$.
We also write $\BN \equiv \{ 1,2,\ldots \}, \BR_+ \equiv [0,\infty)$.
 
We write $X_n \inprob X$ if $X_n$ converges to $X$ in probability, $X_n
\todist X$ if $X_n$ converges to $X$ in distribution, and $X \indist Y$
if $X$ and $Y$ have the same distribution.  We use $\bone$ to denote
indicator functions, so that $\bone(A)=1$ if $A$ is true and
$\bone(A)=0$ otherwise.

 \section{Independent random variables conditioned on a weighted
sum}\label{sect2}

\subsection{The combinatorial setup}

Common to the enumerations (\ref{Nassembly}) through (\ref{Nselection})
 is the form
\eq\label{Ngeneral}
 N(n, \b{a})  \equiv   | \{  \b{C} = \b{a} \} | = 
 \bone(a_1+2 a_2 + \cdots + n a_n =n) 
\  f(n)  \ \prod_1^n g_i(a_i),
\en
with $f(n)=n!$ for assemblies, and $f(n) \equiv 1$ for multisets and
selections. 
To see that (\ref{Ngeneral}) involves independent random variables
conditioned on a weighted sum, view the right hand side as a product of
three factors.
First, the  indicator function, which depends
on both $n$ and $\b{a}$, corresponds to conditioning on the value of a
weighted sum. Second, the factor $f(n)$ does not depend 
on $\b{a}$, and hence disappears from
conditional probabilities.  The product form of the third factor 
corresponds to $n$ mutually independent, but not identically distributed,
random variables. 

The distribution of a random assembly, multiset, or selection $\b{C}(n)$ 
given in (\ref{combdist.old}) can now be expressed in the following form. 
For $\b{a} \in
\B{Z}_+^n$,
\eq\label{combdist}
\B{P}(\b{C}(n)=\b{a})  =
 \bone(a_1+2 a_2 + \cdots + n a_n =n) 
\  \frac{f(n)}{p(n)}  \ \prod_1^n g_i(a_i).
\en

Given functions $g_1,g_2,\ldots$ from $\B{Z}_+ $ to $\B{R}_+$, and a
constant $x>0$, let $Z_1,Z_2, \ldots$ be independent nonnegative
integer valued 
random variables
with distributions given by
\eq\label{not used 1}
\B{P}(Z_i=k) = c_i \ g_i(k) \ x^{ik}, \quad k=0,1,2,\ldots \ .
\en
The above definition, in which $c_i \equiv c_i(x)$ is the normalizing
constant, makes sense if and only if the value of $x$ and the functions
$g_i$ are such that
\eq\label{normconstant}
c_i \ \equiv \ \left( \sum_{k \ge 0} g_i(k) x^{ik} \right)^{-1}
\ \ \in (0,\infty).
\en

For assemblies, $g_i(k) =(m_i/i!)^k/k!$, so that (\ref{normconstant}) is
satisfied for all $x>0$. Defining $\lambda_i \equiv m_i \ x^i/i!$,
we see that $c_i = \exp(-\lambda_i)$ and $Z_i$ is Poisson with mean and
variance 
\eq\label{Zassembly}
\B{E} Z_i = \mbox{var}(Z_i) = \lambda_i \equiv \frac{m_i x^i}{i!}.
\en
For multisets, $g_i(k) = {m_i + k-1 \choose k}$, so the summability
condition (\ref{normconstant}) is satisfied if and only if $x < 1$.  For
$x \in (0,1)$, we have $c_i = (1-x^i)^{m_i}$ and $Z_i$ has the negative
binomial distribution with parameters $m_i$ and $x^i$ given by
\[
\B{P}(Z_i =k) = {m_i +k-1 \choose k} (1-x^i)^{m_i} \ x^{ik}, \ \ \
k=0,1,2,\ldots,
\]
with mean and variance 
\eq\label{Zmultiset}
\B{E} Z_i = \frac{m_i x^i}{1-x^i}, \ \ \ \ \mbox{var}(Z_i) = 
\frac{m_i x^i}{(1-x^i)^2}.
\en
In the special case $m_i=1$, this is just the geometric distribution,
and in general $Z_i$ is the sum of $m_i$ independent random variables 
each with the geometric distribution  $\B{P}(Y=k)=(1-x^i)x^{ik}$ for $k \ge
0$.

For selections, $g_i(k)={m_i \choose k}$, which is zero for $k > m_i$,
so that
(\ref{normconstant}) is satisfied for all $x>0.$  We see that 
 $c_i= (1+x^i)^{-m_i}$, by writing
\[
\B{P}(Z_i=k)= c_i {m_i \choose k} x^{ik} ={m_i \choose k} \left(
\frac{x^i}{1+x^i} \right)^k \left( \frac{1}{1+x^i} \right)^{m_i-k}.
\]
Thus, with $p_i = x^i/(1+x^i)$, the distribution of $Z_i$ is 
binomial with parameters $m_i$ and $p_i$,  with mean and variance
\eq\label{Zselection}
\B{E} Z_i = m_i \ p_i= \frac{m_i x^i}{1+x^i},
 \ \ \ \ \mbox{var}(Z_i)=m_i \ p_i \ (1-p_i)= \frac {m_i x^i}{(1+x^i)^2}.
\en

\subsection{Conditioning on  weighted sums in general}
In order to give a proof of (\ref{equaldist}) which will also serve in
Section \ref{sect6} on process refinements, and Section \ref{sect8} on
large deviations, 
we generalize to a situation
that handles weighted sums with an arbitrary finite index set.  We
assume that $I$ is a finite set, and for each $\alpha \in I,
 g_\alpha : \B{Z}_+ \rightarrow \B{R}_+$ is given. 
We assume that $w$ is a given weight function with values in
 $\B R$ or more generally $\B{R}^d$, so that for $\alpha \in I,
w(\alpha)$ is the weight of $\alpha$.  For the combinatorial examples
in  Section
\ref{sect1}, we had $I=\{1,2,\ldots,n \}$, and a one--dimensional space
of weights, with $w(i)=i$.  
For $\b{a} \in \B{Z}_+^I$ with coordinates $a_\alpha \equiv a(\alpha),$
we use vector dot product notation for the weighted sum
\[
\b{w} \cdot \b{a} \equiv \sum_{\alpha \in I} a(\alpha) w(\alpha) .
\]

Furthermore, we assume that we are given a
target value $t$ such that there exists a normalizing
constant $f(I,t)$ so that the formula
\eq\label{gencombdist}
 \B{P}(\b{C}_I=\b{a}) = \bone(\b{w}\cdot \b{a} = t) f(I,t) \prod_{\alpha
\in I} g_{\alpha}(a_{\alpha}), \quad \b{a} \in \B{Z}_+^I 
\en
defines a  probability distribution for a stochastic process $\b{C}_I$ with
values in  $\B{Z}_+^I$.   
The distribution of $\b{C}(n)$ given by (\ref{combdist}) is a special
case
of (\ref{gencombdist}) with $t=n$ and $f(I,t)= f(n)/p(n).$ 

Assume that for some value $x>0$ there exist normalizing constants
 $c_\alpha \equiv c_\alpha(x) 
\in (0,\infty)$, such that for each $\alpha \in I$,
\begin{equation}
 \ \B{P}(Z_\alpha = k) = c_\alpha(x) g_\alpha (k) x^{w(\alpha)k},
   \ k=0,1,2,\ldots 
\label{genZdist}
\end{equation}
defines a probability distribution on $\B{Z}_+$.  In case $d >1$, so
that $w(\alpha)=(w_1(\alpha),\ldots,w_d(\alpha))$, we take $x \equiv
(x_1,\ldots,x_d) \in (0,\infty)^d$, and $x^{w(\alpha)k}$ denotes the
product $x_1^{w_1(\alpha)k}\cdots x_d^{w_d(\alpha)k}$.
 Define the weighted sum
$T$
by
\eq
\label{Tdef}
T \equiv T_I \equiv \sum_{\alpha \in I} w(\alpha) Z_\alpha.
\en

It should now be clear that the following is a generalization of
(\ref{equaldist}).
\begin{theorem}\label{genequaldist}
Let $\b{Z}_I \equiv (Z_\alpha)_{\alpha \in I}$ have independent 
coordinates $Z_\alpha$
with distribution given by (\ref{genZdist}), and let $\b{C}_I$ have the 
distribution given by (\ref{gencombdist}). Then
\eq
\label{equation genequaldist}
\b{C}_I \indist \left( \b{Z}_I \left| \right. T=t \right),
\en
and hence for any $B \subset I$, the processes restricted to indices in
$B$ satisfy
\eq
\label{equation genequaldist1}
\b{C}_B \indist \left( \b{Z}_B \left| \right. T=t \right).
\en
Furthermore, the normalizing constants and the conditioning probability are 
related by 
\eq
\label{probT=t}
\B{P}(T=t) = f(I,t)^{-1} \  x^t \prod_{\alpha \in I} c_\alpha(x).
\en
\end{theorem}

\no{\bf Remark:} 
The distribution of $\b Z_I$, and hence that of $T \equiv \b
w \cdot \b Z_I$, depends on $x$, so the left side $\B P(T=t)$ of
(\ref{probT=t}) is a function of $x$.

\proof The distribution of $\b{Z}_I$ is given by
\[
 \B{P}(\b{Z}_I=\b{a})= 
 \prod_{\alpha \in I} \left( c_\alpha g_\alpha (a_\alpha)
 x^{w(\alpha)a(\alpha)} \right) =
x^{\b{w} \cdot \b{a}} \ \ \prod_{\alpha \in I} c_\alpha \ \
 \prod_{\alpha \in I} g_\alpha (a_\alpha),
\] 
for  $\b{a} \in \B{Z}_+^I$, so that
if $\b{w} \cdot \b{a} = t$  then 
\eq\label{new1}
\B{P}(\b{Z}_I=\b{a}) = x^t \ \prod c_\alpha \ 
f(I,t)^{-1} \  \B{P}(\b{C}_I=\b{a}).
\en
The conditional distribution of  $\b{Z}_I$ given $\{ T=t \}$ is 
given by
\begin{eqnarray}
\B{P}(\b{Z}_I=\b{a} | T=t) & = &
 \frac{\bone(t= \b w \cdot \b a ) \B{P}(\b{Z}_I=\b{a})}  
   { \B{P}(T=t)}    \nonumber \\
&&\nonumber \\ & = &
 \frac{x^t \ (\prod c_\alpha) f(I,t)^{-1} \B{P}(\b{C}_I=\b{a})}  
   { \B{P}(T=t)}
\label{step1} \\
&&\nonumber \\& = &
 \frac{\ \ \ \ x^t \ (\prod c_\alpha) f(I,t)^{-1} \B{P}(\b{C}_I=\b{a})}  
   { \sum_{\b{b} \in \B{Z}_+^I}x^t \ (\prod c_\alpha) \ 
f(I,t)^{-1}
 \ \B{P}(\b{C}_I=\b{b})}
\nonumber\\
&&\nonumber \\& = &
\frac{ \B{P}(\b{C}_I=\b{a})}  
   { \sum_{\b{b} } \B{P}(\b{C}_I=\b{b})} \nonumber \\
&&\nonumber \\
& = & \B{P}(\b{C}_I=\b{a}), \ \b{a} \in \B{Z}_+^I. 
\label{step2}
\end{eqnarray}
The equality between (\ref{step1}) and (\ref{step2}), for any $\b{a}$ for
which $\B{P}(\b{C}_I=\b{a})>0$, establishes (\ref{probT=t}).\hfill \qed

For the combinatorial objects in Section {\ref{sect1}, $I=\{
1,2,\ldots,n \}$, and $w(i)=i$.  For this case $T$ reduces to
\eq
    T \equiv T_n \equiv Z_1 + 2 Z_2 + \cdots + n Z_n.
\label{deftn}
\en
In the  case of assemblies, corresponding to (\ref{Nassembly}) and
(\ref{Zassembly}),  the distribution of $Z_i$ is Poisson $(\lambda_i)$, 
and (\ref{probT=t}) reduces to
\eq
\label{assprobT=t}
\B P (T_n=n) \  = \ \frac{p(n)}{n!} \ x^n \ \exp(-\lambda_1-\cdots
-\lambda_n),
\en
where  $\lambda_i = m_i x^i /i!$ and $ x>0$.
In the case of multisets, corresponding to (\ref{Nmultiset})
and (\ref{Zmultiset}),  $Z_i$ is distributed like 
the sum of $m_i$ independent
geometric $(x^i)$ random variables, and (\ref{probT=t}) reduces to
\eq
\label{multiprobT=t}
\B P (T_n=n)= \ p(n) \ x^n \ \prod_1^n (1-x^i)^{m_i},
\en
for $0<x<1$.
In the case of selections, corresponding to (\ref{Nselection}) and
(\ref{Zselection}), the distribution of  $Z_i$
is binomial $(m_i, x^i/(1+x^i))$, so that (\ref{probT=t}) reduces to
\eq
\label{selectprobT=t}
\B P (T_n=n) = p(n) \ x^n \ \prod_1^n (1+x^i)^{-m_i},
\en
for $x>0$.

 
\section{Total variation distance}\label{sect3}

A useful way to establish that the independent process $\b Z_n \equiv
(Z_1,Z_2,\ldots,Z_n)$ is a good approximation for the dependent
combinatorial process $\b C(n)$ is to focus on a subset $B$ of the
possible component sizes, and give an upper bound on the total variation
distance between the two processes, both restricted to $B$.  Theorem
\ref{tvthm} below shows how this total variation distance for
these two processes reduces to the total variation distance between 
two one--dimensional random variables.  

Here is a quick review of the relevant features of total variation
distance. For two random elements $X$ and $Y$ of a finite or countable
space $S$,
the total variation distance between $X$ and $Y$ is defined by
\[
d_{TV}(X,Y) = \frac{1}{2} \sum_{s \in S} |\B P(X=s) - \B P(Y=s)|.
\]
Properly speaking this should be referred to as the distance between the
distribution ${\mathcal L}(X)$ of $X$ and the distribution ${\mathcal L}(Y)$ of
$Y$, written for example as $d_{TV}({\mathcal L}(X),{\mathcal L}(Y))$.
Throughout this paper we use the simpler notation, except in Section
\ref{sect8} which involves changes of measure.

Many authors, following the tradition of analysis of signed measures,
 omit the factor of $1/2$.  Using the factor of
$1/2$, we have that $d_{TV}(X,Y) \in [0,1]$, and furthermore, $d_{TV}$ is
identical to the Prohorov metric, providing the underlying metric on $S$
assigns distance $\ge 1$ between any two distinct points. 
In particular,
a sequence of random elements $X_n$ in a discrete space $S$ converges in
distribution to $X$ if and only if $d_{TV}(X_n,X) \rightarrow 0$. 
 
Another characterization of total variation distance
 is
\[
d_{TV}(X,Y) = \max_{A \subset S}\left( \B P(X \in A) - \B P(Y \in A)
\right),
\]
and in the discrete case, 
a necessary and sufficient condition that the maximum be achieved by $A$
is that $\{s: \B P(X=s) > \B
P(Y=s) \}  \subset A \subset \{s: \B P (X=s) \ge \B P(Y=s) \}$.

The most intuitive description of total variation distance is in terms
of coupling.  A
``coupling'' of $X$ and $Y$ is a probability measure on $S^2$ whose
first and second marginals are the distributions of $X$ and $Y$
respectively.  Less formally, a coupling of $X$ and $Y$ is a 
recipe for constructing $X$ and $Y$ simultaneously on the same
probability space, subject only to having given marginal distributions for
$X$ and for $Y$. In terms of all possible coupling measures on $S^2$,
\eq\label{tvcouple}
d_{TV}(X,Y) = \min_{couplings} \B P (X \ne Y).
\en
The minimum above is achieved, but
 in general there is not a unique optimal coupling. In fact
a discrete coupling achieves $\B P (X \ne Y) = d_{TV}(X,Y)$, if and only if, for
all $s \in S, \ \B P(X=Y=s) = \min ( \B P(X=s),\B P (Y=s))$. 
Intuitively, if  $d_{TV}(X,Y)$  is small, then $X$ and $Y$ are nearly
indistinguishable from a single observation; formally, for any
statistical test to decide whether $X$ or $Y$ is being observed, the sum
of the type I and type II errors is at least $1-$  $d_{TV}(X,Y)$ .

Upper bounds on the total variation distance between a combinatorial
process and a simpler process are useful because these upper bounds are
inherited by functionals of the processes.  If  $h:S
\rightarrow T$ is a deterministic map between countable spaces,
 and $X$ and $Y$ are random
elements of $S$, so that $h(X)$ and
$h(Y)$ are random elements of $T$, then
\eq
d_{TV}(h(X),h(Y)) \ \leq \ d_{TV}(X,Y).
\label{tvinequality}
\en
Theorem \ref{tvthm} below, and its refinement, Theorem \ref{tvrefine} in
Section \ref{sect6}, both describe combinatorially interesting cases in
which  equality holds in (\ref{tvinequality}).  It is natural to
ask when, in general, such equality holds. The following elementary
theorem provides an answer. 

\begin{theorem}\label{tv=thm}
In the discrete case, equality holds in (\ref{tvinequality}) if and only
if the sign of $\B P(X=s)- \B P(Y=s)$ depends only  on $h(s)$, in the
non-strict sense that  
 $\forall a,b \in S$,
\[ h(a)=h(b) \mbox{ implies } 
\left(\B P(X=a) - \B P(Y=a)\right)\ \left(\B P(X=b)-\B P(Y=b)\right) \ \ \ge 0.
\]
\end{theorem}
\proof
Consider the proof of (\ref{tvinequality}), namely
\begin{eqnarray}
2 d_{TV}(h(X),h(Y))& =& \sum_{r \in T} \left| \B P(h(X)=r)-\B P(h(Y)=r)
\right| \nonumber \\
&&\nonumber \\ 
&=& \sum_r \left| \sum_{a \in S: h(a)=r} (\B P(X=a)-\B P(Y=a))
\right| \label{cancelstep} \\
&& \nonumber 
\\ &\leq & \sum_r  \sum_{a: h(a)=r}\left| \B P(X=a)-\B P(Y=a)
\right| \label{canstep} \\
&& \nonumber\\ 
&=& 2 d_{TV}(X,Y). \nonumber  
\end{eqnarray}
Since the inequality in (\ref{canstep}) holds term by term in the outer
sums,
equality holds overall if and only if equality holds for each $r$. This
 in turn is equivalent to the condition that   for each $r$,
 there are no terms of opposite sign in the inner sum
in (\ref{cancelstep}).\hfill \qed

Diaconis and Pitman (1986) view ``sufficiency'' as a
key concept.  In the context above, $h: S \rightarrow T$ is a sufficient
statistic for discriminating between the distributions of $X$ and $Y$ in
$S$, if the likelihood ratio depends only on $h$; i.e. if there is a
function $f: T \rightarrow \B R$ such that for all $s \in S, \ \B P(X=s)
= f(h(s)) \ \B P(Y=s)$. Taking a sufficient statistic preserves total
variation distance, as observed by
 Stam (1978). This is also a special case of Theorem 
\ref{tv=thm}, in
which a product is nonnegative because it is a square: $ (\B P(X=a)-\B
P(Y=a))(\B P(X=b)-\B P(Y=b)) \ = \ (f(h(a))-1)(f(h(b))-1)\B P(Y=a) \B
P(Y=b) \geq 0$ whenever $h(a)=h(b)$.
 
\begin{theorem}\label{tvthm}
Let $I$ be a finite set, and for $\alpha \in I$, let $C_\alpha$ and
$Z_\alpha$ be $\B Z_+$ valued random variables, such that the
$Z_\alpha$ are mutually independent.  Let $\b w = (w(\alpha))_{\alpha \in
I}$ be a deterministic weight function on $I$ with values in some linear
space, let $T = \sum_{\alpha \in
I} w(\alpha) Z_\alpha$, and let $t$ be such that
$\ \B P(T=t) >0$.  For $B \subset I$, 
we use the notation $\b C_B \equiv (C_\alpha)_{\alpha \in
B}$ and $\b Z_B \equiv (Z_\alpha)_{\alpha \in B}$ for
random elements of $\B Z_+^B$.  Define 
\[
R \equiv R_B \equiv \sum_{\alpha \in B} w(\alpha) Z_\alpha,\ \  
S \equiv S_B \equiv \sum_{\alpha \in I-B} w(\alpha) Z_\alpha,
\]
so that $T = R+S$ and $R$ and $S$ are independent.
If 
\eq
\b C_I \indist (\b Z_I | T=t),
\label{tvhyp}
\en 
 then
\eq\label{tvhyp1}
d_{TV}(\b C_B, \b Z_B) \ = d_{TV}( (R_B|T=t), R_B).
\en
\end{theorem}
\proof
We present two proofs, since it is instructive to contrast them.
\def\rofa{\b w \cdot \b a}
Note that not only are $R$ and $S$ independent, but also that $R$ is a
function of $\b Z_B$, and $\b Z_B$ and $S$ are independent. 
 For $\b a \in \B Z_+^B$, write $\rofa \equiv \sum_{\alpha \in B}
w(\alpha) a(\alpha)$. 
\def\fh{\frac{1}{2}}
\begin{eqnarray*}
d_{TV}(\b C_B, \b Z_B) & = & \fh \sum_{\b a \in \B Z_+^B}
 \left| \B P( \b Z_B = \b a \left| T=t \right. ) - 
\B P( \b Z_B = \b a) \right| \\
& = &
   \fh \sum_r \ \ \ \sum_{\b a: \rofa = r}
 \left| \frac{\B P( \b Z_B = \b a, r+S=t)}{\B P(T=t)}
   - \B P( \b Z_B =\b a) \right|
\\ && \\ & = &
   \fh \sum_r \ \ \ \sum_{\b a: \rofa = r}
 \left| \frac{\B P( \b Z_B = \b a) \B P( r+S=t)}{\B P(T=t)}
   - \B P( \b Z_B =\b a) \right|
\\ && \\ & = &
   \fh \sum_r 
 \left| \frac{\B P( R=r) \B P( r+S=t)}{\B P(T=t)}
   - \B P( R=r ) \right|
\\ && \\ & = &
   \fh \sum_r \left| \frac{\B P(R=r,r+S=t)}{\B P(T=t)}
   - \B P(R=r) \right|
\\ && \\ &=&
   \fh \sum_r \left| \B P(R=r|T=t) - \B P(R=r) \right|
\\ && \\ &=& d_{TV}( (R|T=t), R).
\end{eqnarray*}

Here is a second proof of  Theorem \ref{tvthm}, viewed 
as a corollary of Theorem \ref{tv=thm}, with the functional $h$ on $\B
Z_+^B$ defined by $h(\b a) = \rofa$.  We need only observe that $h$ is a
sufficient statistic since $\B P(\b Z_B=\b a|T=t)= \BP(\b Z_B=\b a)
\BP(S=t-h(\b a))/\BP(T=t) \ $.   \hfill \qed

For the sake of calculations of total variation distance between a
combinatorial process and its independent process approximation, the
most useful form for the conclusion of Theorem \ref{tvthm} is
\begin{eqnarray}
d_{TV}(\b C_B, \b Z_B) & = &
   \frac{1}{2} \sum_r 
 \left| \frac{\B P( R=r) \B P( r+S=t)}{\B P(T=t)}
   - \B P( R=r ) \right| \nonumber
\\ && \nonumber \\ & = &
   \frac{1}{2} \sum_r \B P( R=r) 
 \left| \frac{\B P( S=t-r)}{\B P(T=t)}
   - 1 \right|. \label{tvcalc}
\end{eqnarray}
In the  usual combinatorial case, where $t=n$ and
$T=Z_1+2Z_2+\cdots+nZ_n$, this gives 
\eq
d_{TV}(\b C_B, \b Z_B)=
 \frac{1}{2} \B P(R>n) +\frac{1}{2} \sum_{r=0}^n \B P( R=r) 
 \left| \frac{\B P( S=n-r)}{\B P(T=n)}
   - 1 \right|. \label{tvcalc1}
 \en

There are two elementary observations that point to strategies for giving
upper bounds on total variation distance.  First, for discrete random
elements we have in general
\begin{eqnarray*}
d_{TV}(X,Y) & \equiv & \frac{1}{2} \sum_{s \in S} |\B P(X=s) - \B P(Y=s)|
\\ && \\
&=& \sum_{s \in S} \left( \B P(X=s) - \B P(Y=s) \right)^+
\\ && \\
&=& \sum_{s \in S} \left( \B P(X=s) - \B P(Y=s) \right)^-,
\end{eqnarray*}
where the notation for positive and negative parts is such that, for
real $x$, $x= x^+ - x^-,$ and $|x|=x^+ + x^-$.  In the context of 
(\ref{tvcalc}) this is
useful in the following form. Let $A \subseteq I$. Then
\begin{eqnarray}
d_{TV}(\b C_B, \b Z_B) & = &
   \sum_r \B P( R=r) 
 \left(1 - \frac{\B P( S=t-r)}{\B P(T=t)}
    \right)^+ \nonumber \\ && \nonumber \\
& \leq & \B P(R \not \in A) \ \ + \sup_{r \in A} 
      \left(1 - \frac{\B P( S=t-r)}{\B P(T=t)}
    \right)^+. \label{truncstrategy} 
\end{eqnarray}
Specializing to the case where the weighted sum $R$ is real
valued, and $A = \{0,1,2,\ldots,k\}$, the truncation level $k$ 
is chosen much larger than $\B E R$,
so that large deviation theory can be used to bound $\B P(R > k)$, but
not too large, so that $\B P( S=t-r)/ \B P(T=t)$ can be controlled to
show it is close to one.  

The second elementary observation, which is proved and exploited in
\TVone, is that the denominator in (\ref{tvcalc}) can be replaced by any
constant $c>0$, at the price of at most a factor of 2, in the sense that
for independent $R$ and $S$ such that $\B P(R+S=t) >0$,
\[
\frac{1}{2} \sum_r \B P(R=r) \left|\frac{\B P(S=t-r)}{\B P(R+S=t)} -1
\right| 
\leq  
 \sum_r \B P(R=r) \left|\frac{\B P(S=t-r)}{c} -1 \right| .
\]
By using this, for example with $c=\B P(S=t)$, 
giving an upper bound on the total variation distance for
combinatorial process approximations is reduced to showing that the
density of $S$ is relatively constant.

Lower bounds for variation distance are often more difficult to obtain,
but it is worth noting that in the combinatorial setup, 
since $\{R_B>n\} \subseteq \{\bC_B \neq
\bZ_B\}$, we have, without the factor $\fh$ suggested by
(\ref{tvcalc1}),
\eq\label{tvlower}
d_{TV}(\b C_B,\b Z_B) \ge \BP(R_B >n).
\en

  \section{Heuristics for useful approximation}\label{sect4}
 
 Recall first that for $B \subset [n]$, we have  
 $\b C_B \indist (\b Z_B | T_n = n)$.
 If $d_{TV}(\b C_B, \b Z_B)$ is small, the approximation of $\b C_B$ by
 $\b Z_B$ is useful. Probabilistic intuition suggests that conditioning
 on $T_n = n$ does not change the distribution of $\b Z_B$ by much,
 provided that the event $\{T_n = n\}$ is relatively likely. This in turn
 corresponds to a choice of $x = x(n)$ for which $\BE T_n$ is
 approximately $n$. Let $\sigma_n^2 \equiv {\rm var}(T_n),$
 and  let $\sigma_B^2 = {\rm var}(R_B)$.
 Intuition then suggests that {\em if}
 \eq\label{metathm1}
 \frac{n - \BE(T_n)}{\sigma_n} \mbox{\it\ is not large}
\en
{\em and}
\eq\label{metathm1a}
  \frac{\BE R_B}{\sigma_n}
 \mbox{\it\ and }\frac{\sigma_B}{\sigma_n} \mbox{\it\ are small}
 \en
 {\em then} $d_{TV}(\b C_B, \b Z_B)$ {\em is small.}

 While our main focus is on the appropriate choice of $x$, we also
 discuss below the appropriate choice of $B$ for examples including
 permutations, mappings, graphs, partitions of sets, and partitions of
 integers.
 
 There is an important qualitative distinction between cases in which the
 appropriate $x$ is constant, and those in which $x$ varies with $n$.
 If $x$ does not depend on $n$, then a single independent process $\b Z =
 (Z_1,Z_2,\ldots)$ may be used to approximate $\b C(n) \equiv
 (C_1(n),\ldots,C_n(n))$, which we identify with
 $(C_1(n),\ldots,C_n(n),0,0,\ldots) \in \BZ_+^\infty$.
 Under the usual product topology on $\BZ_+^\infty,$ we have that 
  $\b C(n) \Rightarrow \b Z$  if, and only if,  for
 every fixed $b$, $\b C_b(n) \equiv
 (C_1(n),\ldots,C_b(n)) \Rightarrow \b Z_b \equiv (Z_1,\ldots,Z_b)$
 as random elements in $\B Z_+^b$.
  Since the metric on $\BZ_+^b$ is discrete, we conclude
 that $\b C_b(n) \Rightarrow \b Z_b$ if, and only if, for each fixed $b$, 
 $d_{TV}(\b C_b(n),\b Z_b) \to 0$. For cases where $x$, and hence $\b Z$, varies
 with $n$, it makes no sense to write 
 $\b C(n) \Rightarrow \b Z$.
 However, it is still useful to be able to
 estimate $d_{TV}(\b C_B(n),\b Z_B(n))$.
 
 We discuss first considerations involved in the choice of $x$ and $B$, and then
 heuristics for predicting the accuracy of approximation.
 
 \subsection{Choosing the free parameter $x$}
 
 It is convenient to discuss the three basic types of combinatorial
 structure separately.
 
 \subsubsection{Assemblies}
 It follows from (\ref{Zassembly}) that
 \eq\label{assem-etn}
 \BE T_n \equiv \sum_{i=1}^n i \BE Z_i = \sum_{i=1}^n \frac{m_i x^i}{(i-
 1)!},
 \en
 while 
 \eq\label{assem-vartn}
 \sigma_n^2 = \sum_{i=1}^n i^2 \BE Z_i = \sum_{i=1}^n \frac{i^2 m_i x^i}{i!}.
 \en
 In the case of permutations, we take $x = 1$ to see that $\BE T_n = n$,
 and $\sigma_n^2 = n(n+1)/2$. In \TVone\ it is proved that $d_{TV}(\b C_B,\b
 Z_B) \rightarrow 0$ as $n\rightarrow \infty,$ with $B=B(n)$, if and
 only if $|B|=o(n)$.  

For the class of assemblies which satisfy the additional condition
 \eq\label{assembly-hyp}
 \frac{m_i}{i!} \sim \frac{\kappa y^i}{i} {\rm\ as\ }i \to \infty,
 \en
 where $y > 0$ and $\kappa > 0$ are constants, we see that
$$
\frac{\BE T_n}{n} \to \left\{
\begin{array}{cl}
  0, & {\rm \ if\ } 0 < x < y^{-1} \\
  \kappa, & {\rm \ if\ } x = y^{-1}\\
 \infty, & {\rm \ if\ } x > y^{-1}.
\end{array}
\right.
$$
 Hence the only fixed $x$ that ensures that $\BE T_n \bothsides n$ is
 $x = y^{-1}$, in which case 
 \eq\label{mean-varass}
 \BE T_n \sim n \kappa, \  \sigma_n \sim n
 \sqrt{\frac{\kappa}{2}}.
 \en 
 
 For the example of random mappings,
 $$
 m_i = e^i (i-1)!\, \BP({\rm Po}(i) < i),
 $$
 where ${\rm Po}(i)$ denotes a Poisson random variable with mean $i$;
see Harris (1960), Stepanov (1969).
 It follows that we must take $x = 1/e$, and, from the Central Limit
 Theorem, $\kappa = 1/2$. 
 In this case $\BE T_n \sim n/2$ and $ \sigma_n \sim n/2$.
 
 For the example of random graphs, with all $2^{{n \choose 2}}$ graphs
 equally likely, the fact that the probability of being connected tends
 to 1 means that the constant vector $(0,0,\ldots,0,1) \in \B Z_+^n$ is a
 good approximation, in total variation distance, to $\b C(n)$.  This is
 a situation in which the equality $\b C(n) \indist (\b Z_n|T_n=n)$
 yields no useful approximation.  With $x$ chosen so that $\B ET_n=n$,
 and $B=\{1,2,\ldots,n-1\}$, we have that $d_{TV}(\b C_B,\b Z_B)
 \rightarrow 0$, but only because both distributions are close to that
of the process that is identically 0 on $\B Z_+^B$.
 
 For partitions of a set, which is discussed further in Section
 \ref{sect5.2} and Section \ref{sect10}, with $x=x(n)$ being the solution of
 $xe^x=n$, and $B=\{1,2,\ldots,b \} \cup \{c,c+1,\ldots,n \}$
 where $b \equiv b(n)$ and $c \equiv c(n)$,
 the heuristic (\ref{metathm1}) suggests that $d_{TV}(\b C_B, \b Z_B)
 \rightarrow 0$ if and only if both $(x-b)/\sqrt{\log n} \rightarrow
 \infty$ and $(c-x)/\sqrt{\log n}
 \rightarrow \infty$.
 For $B$ of the complementary form
 $B=\{b,b+1,\ldots,c \}$ with $b<c$ both within a bounded number of
 $\sqrt{\log n}$ of $x$,  the heuristic suggests that $d_{TV}(\b C_B, \b Z_B)
 \rightarrow 0$ if, and only if, $(c-b)=o(\sqrt{\log n})$. Sachkov
(1974) and Fristedt (1992) have partial results in this area.
 
 \subsubsection{Multisets}
 Using (\ref{Zmultiset}) we see that
 \eq\label{multi-etn}
 \BE T_n  = \sum_{i=1}^n \frac{i m_i x^i}{1-x^i},
 \en
 while 
 \eq\label{multi-vartn}
 \sigma_n^2 = \sum_{i=1}^n \frac{i^2 m_i x^i}{(1-x^i)^2}.
 \en
 If the multiset construction satisfies the additional hypothesis that
 \eq\label{multiset-hyp}
 m_i \sim \frac{\kappa y^i}{i}\ {\rm \ as\ } i \to \infty,
 \en
 where $y > 1$ and $\kappa > 0$ is fixed, a similar analysis shows
 that the only fixed $x$ that ensures that $\BE T_n \bothsides n$ is
 $x = y^{-1}$, in which case the asymptotics for $\BE T_n$ and
$\sigma_n$ are the same as those in (\ref{mean-varass}).
 
 The first example that satisfies the hypothesis in (\ref{multiset-hyp})
 is the multiset in which $p(n) = q^n$ for some integer $q \geq 2$. In this case
 the $m_i$ satisfy
 \eq\label{poly1}
 q^n = \sum_{j|n} j m_j,
 \en
 so that by the M\"obius inversion formula we have
 \eq\label{poly2}
 m_n = \frac{1}{n} \sum_{k|n} \mu(n/k) q^k,
 \en
 where $\mu(\cdot)$ is the M\"obius function, defined by
 \begin{eqnarray*}
 \mu(n) & = & (-1)^k {\rm\ if\ } n {\rm\ is\ the\ product\ of\ } k {\rm 
 \ distinct\ primes}\\
 \mu(n) & = & 0 {\rm\ otherwise}.
 \end{eqnarray*}
 It follows from (\ref{poly1}) that
 $$
 q^i - \frac{q}{q-1} q^{i/2} \leq i m_i \leq q^i,
 $$
 so that (\ref{multiset-hyp}) holds with $\kappa = 1, y = q$. This
 construction arises in the study of necklaces (see Metropolis and Rota
 (1983, 1984) for example), in card shuffling (Diaconis, McGrath and
 Pitman (1994)), and, for $q$ a prime power, in 
 factoring polynomials over $GF(q)$, a finite field of $q$ elements. 
 In this last case $m_i$ is  the number of irreducible monic polynomials over
 $GF(q)$; see Lidl and Niederreiter (1986), for example. 
 
 Another example concerns random mapping patterns. Let $t_n$ denote the number of rooted trees with $n$ unlabelled points, and set $T(x) = \sum_{n=1}^\infty t_n x^n.$ Otter (1948) showed that $T(x)$ has radius of convergence $\rho = 0.3383\ldots$, from which Meir and Moon (1984) established that
 $$
 m_i \sim \frac{\rho^{-i}}{2i}.
 $$
  Hence (\ref{multiset-hyp}) applies with $\kappa = 1/2, y = \rho^{-1}$.
 
 For an example in which $x$ varies with $n$, we consider random
 partitions of the integer $n$. In this case $m_i \equiv 1$. Taking 
 $x = e^{-c/\sqrt{n}}$
 and using (\ref{multi-etn}), we see that
 \begin{eqnarray*}
 n^{-1} \BE T_n & = & \sum_{i=1}^n \frac{n^{-1/2} i \exp(-i
 c/\sqrt{n})}{1 - \exp(-i c / \sqrt{n})} \frac{1}{\sqrt{n}} \\
 & \to & \int_0^{\infty} \frac{y e^{-c y}}{1 - e^{- c y}} dy \\
 & = & \frac{1}{c^2} \int_0^1 \frac{- \log(1-v)}{v} dv \\
 & = & \frac{\pi^2}{6 c^2}.
 \end{eqnarray*}
 Hence to satisfy $\BE T_n \sim n$, we choose $c = \pi/\sqrt{6}$, so that
 \eq\label{choosex-multi}
 x = \exp(-\pi/\sqrt{6 n}).
 \en
 From (\ref{multi-vartn}), it follows by a similar calculation that
 \begin{eqnarray*}
 n^{-3/2}\, \sigma_n^2 & \to & \int_0^{\infty} \frac{y^2 e^{-c y}}
 {(1 - e^{- c y})^2} dy \\
 & = & \frac{1}{c^3} \int_0^1 \left(\frac{- \log(1-v)}{v}\right)^2 dv \\
 & = & \frac{2}{c},
 \end{eqnarray*}
 so that
 \eq\label{vartn-multipart}
 \sigma_n^2 \sim \frac{2 \sqrt{6}}{\pi} n^{3/2}.
 \en
 For sets of the form $B=\{1,2,\ldots,b \} \cup \{c,c+1,\ldots,n \}$
 where $0 \leq b \equiv b(n)$ and $c \equiv c(n) \leq n$,
 the heuristic in (\ref{metathm1}) and (\ref{metathm1a}) suggests 
that $d_{TV}(\b C_B, \b Z_B)
 \rightarrow 0$ if, and only if, both $b=o(\sqrt{n})$ and $c/\sqrt{n}
 \rightarrow \infty$. For $B$ of the complementary form
 $B=\{b,b+1,\ldots,c \}$ with $b<c$ both of the order of $\sqrt{n}$,
   the heuristic suggests that $d_{TV}(\b C_B, \b Z_B)
 \rightarrow 0$ if, and only if, $(c-b)=o(\sqrt{n})$. See Fristedt
(1993) and Goh and Schmutz (1993) for related results.
 
 \subsubsection{Selections}
 In this case, it follows from (\ref{Zselection}) that
 \eq\label{select-etn}
 \BE T_n  = \sum_{i=1}^n \frac{i m_i x^i}{1+x^i},
 \en
 while 
 \eq\label{select-vartn}
 \sigma_n^2 = \sum_{i=1}^n \frac{i^2 m_i x^i}{(1+x^i)^2}.
 \en
 If the selection construction satisfies the additional hypothesis
 (\ref{multiset-hyp}),
 then, just as for the assembly and multiset constructions, we take 
 $x = y^{-1}$, and (\ref{mean-varass}) holds once more.
 As an example, for square--free factorizations of polynomials over
 a finite field with $q$ elements, we have $y = q, \kappa = 1, x =
 q^{-1}$.
 
 For an example in which $x$ varies with $n$, we consider once more random
 partitions of the integer $n$ with all parts distinct,
 which is the selection construction with $m_i \equiv 1$. Taking 
 $x = e^{-d/\sqrt{n}}$,
 and using (\ref{select-etn}) we see that
 \begin{eqnarray*}
 n^{-1} \BE T_n & = & \sum_{i=1}^n \frac{n^{-1/2} i \exp(-i
 d/\sqrt{n})}{1 + \exp(-i d / \sqrt{n})} \frac{1}{\sqrt{n}} \\
 & \to & \int_0^{\infty} \frac{y e^{-d y}}{1 + e^{- d y}} dy \\
 & = & \frac{1}{d^2} \int_0^1 \frac{- \log v}{1+v} dv \\
 & = & \frac{\pi^2}{12 d^2}.
 \end{eqnarray*}
 Hence to satisfy $\BE T_n \sim n$, we pick $d = \pi/\sqrt{12}$, so that
 \eq\label{choosex-select}
 x = \exp(-\pi/ \sqrt{12 n}).
 \en
 From (\ref{multi-vartn}), it follows by a similar calculation that
 \begin{eqnarray*}
 n^{-3/2}\, \sigma_n^2 & \to & \int_0^{\infty} \frac{y^2 e^{-d y}}
 {(1 + e^{- d y})^2} dy \\
 & = & \frac{1}{d^3} \int_0^1 \left(\frac{- \log v}{1+v}\right)^2 dv \\
 & = & \frac{2}{d}.
 \end{eqnarray*}
 For the choice of $x$ in (\ref{choosex-select}), we see that
 \eq\label{vartn-selectpart}
 \sigma_n^2 \sim \frac{4 \sqrt{3}}{\pi} n^{3/2}.
 \en
  
 To see how easy the heuristic for choosing $x$ can be, consider
 partitions of the integer $n$ with all parts distinct and odd.  Compared
 to the above calculations, we are simply leaving out every other term,
 so that $n^{-1} \B ET_n \rightarrow \pi^2/(24d^2)$, and we prescribe
 using $x= \exp(-\pi/ \sqrt{24n})$.  As with unrestricted partitions, using
 the appropriate $x$ for either partitions with distinct parts or
 partitions with distinct odd parts, we believe that the unconditioned
 process $\b Z_B$ is a good approximation for the combinatorial process
 $\b C_B$, in the total variation sense, if and only if $b/\sqrt{n}$ is
 small and $c/\sqrt{n}$ is large, for 
  $B=\{1,2,\ldots,b \} \cup \{c,c+1,\ldots,n \}$.
 For $B$ of the complementary form
 $B=\{b,b+1,\ldots,c \}$ with $b<c$ both of the order of $\sqrt{n}$, 
  the heuristic suggests that $d_{TV}(\b C_B, \b Z_B)
 $ is small if, and only if, $(c-b)$ is small relative to $\sqrt{n}$.
 
 \subsection{A quantitative heuristic}\label{sect4.2}

 In several examples, the $\BZ_+$-valued random variables $T_n$,
 appropriately centered and rescaled, converge in distribution to a continuous 
 limit $X$ having a density $f$ on $\BR$. For illustration, we 
describe the important class of cases in which 
 \eq\label{rescale0}
 \frac{T_n}{n} \Rightarrow X.
 \en
A local limit heuristic suggests the approximation
\eq\label{local1}
\BP(T_n = n)  \doteq \frac{f(1)}{n},
\en
where the sense of the approximation $\doteq$ is deliberately vague.
Assuming that $B$ is small, so that $R/n \inprob 0$, we also have $S/n
\Rightarrow X$.
For $0 \leq k \leq n$, the local limit heuristic gives
$$
\BP(S = n-k)  \doteq \frac{1}{n}\,f\left(1 - \frac{k}{n}\right),
$$
and a Taylor expansion further simplifies this to 
\eq\label{local2}
\BP(S = n-k) \doteq \frac{1}{n} \left(f(1) - \frac{k}{n} f^\prime(1-
)\right).
\en

Using these approximations in the total variation  formula (\ref{tvcalc}) gives
 \begin{eqnarray*}
 d_{TV}(\b C_B,\b Z_B) & = & \frac{1}{2} \sum_{k = 0}^n \BP(R=k) \left|
 1 - \frac{\BP(S=n-k)}{\BP(T_n = n)}\right| + \frac{1}{2}\BP(R > n)\\
 &&  \\
 & \doteq & \frac{1}{2} \sum_{k \geq 0} \BP(R=k) \left| 1 - \frac{
 n^{-1}  (f(1) - n^{-1} k f^{\prime}(1-)}
 {n^{-1} f(1)}\right| \\
 &&  \\
 & = & \frac{1}{2} \frac{|f^{\prime}(1-)|}{f(1)} 
 \frac{\BE |R|}{n}. 
 \end{eqnarray*}
 However, this approximation ignores the essential feature that
$d_{TV}(\mu,\nu) = \frac{1}{2} |\mu - \nu|$, where the signed measure
 $\mu -\nu$ has net mass zero.
Thus, even though $f(1)/n$ is the natural approximation for $\BP(T_n =
n)$, it is important to use a more complicated heuristic in which the 
approximation for $T$ is the
convolution of the distribution of $R$ and our approximation for the
distribution of $S$. Thus
\begin{eqnarray}
\BP(T = n) & = & \sum_{k =0}^n \BP(R=k) \BP(S = n-k) \nonumber\\
&&\nonumber\\
& \doteq & \sum_{k \geq 0} \BP(R = k) \frac{1}{n} \left(f(1) - \frac{k}{n}
f^\prime(1-) \right) \nonumber \\
&&\nonumber \\
& = & \frac{1}{n} \left( f(1) - \frac{\BE R}{n} f^\prime(1-)\right).
\label{pt=napprox}
\end{eqnarray}
Using this approximation,
 \begin{eqnarray}
 d_{TV}(\b C_B,\b Z_B) & = & \frac{1}{2} \sum_{k = 0}^n \BP(R=k) \left|
 1 - \frac{\BP(S=n-k)}{\BP(T_n = n)}\right| + \frac{1}{2} \BP(R > n)\nonumber\\
 && \nonumber \\
 & \doteq & \frac{1}{2} \sum_{k \geq 0} \BP(R=k) \left| 1 - \frac{
 n^{-1}  (f(1) - n^{-1} k f^{\prime}(1-))}
 {n^{-1} (f(1) - n^{-1} \BE R f^\prime(1-))}\right| \nonumber \\
 && \nonumber \\
& = & \frac{1}{2} \sum_{k \geq 0} \BP(R=k) \left| \frac{n^{-1} (k - \BE R)
f^\prime(1-)}{f(1) - n^{-1} \BE R f^\prime(1-)}\right| \nonumber \\
&& \nonumber \\
& = & \frac{1}{2n} |f^\prime(1-)| \  \BE | R - \BE R|
\ |f(1) - n^{-1} \BE R f^\prime(1-)|^{-1} \nonumber \\
&& \nonumber \\
& \doteq & \frac{1}{2n} \frac{|f^\prime(1-)|}{f(1)} \BE | R - \BE R|.
  \label{approx2}
 \end{eqnarray}

As a plausibility check, we note that the alternative approximation
using $\BP(T_n = n) \doteq \frac{1}{n} f(1)$ and $S \doteq T - \BE R$,
so that $\BP(S = n-k) \doteq \BP(T = n + \BE R - k) \doteq \frac{1}{n}
f(1 - \frac{k-\BE R}{n})$, also satisfies the convolutional property, and
leads to the same first order result as (\ref{approx2}).

 One possible specific interpretation  of the approximation in (\ref{approx2})
 would be the following pair of statements, giving a decay rate for
$d_{TV}$, for fixed $B$, as $n \to
\infty$.
 
 {\it If $T_n/n \Rightarrow X$, and $X$ has density $f$ with
$f^\prime(1-) \neq 0$, then}
 \eq\label{approx3}
 d_{TV}(\b C_B,\b Z_B) \sim \frac{1}{2} \frac{|f^{\prime}(1-)|}{f(1)}
 \frac{\BE|R - \BE R|}{n},
 \en
 
 {\it If $T_n/n \Rightarrow X$, and $X$ has density $f$ with
$f^\prime(1-) = 0$, then}
 \eq\label{approx4}
 d_{TV}(\b C_B,\b Z_B) = o\left(\frac{1}{n}\right).
 \en
 
For the more general case in which there are constants $s_n$ such that
$$
\frac{T_n - n}{s_n} \Rightarrow X
$$
where $X$ has density $f$, these statements are to be replaced by
\eq\label{approx3gen}
 d_{TV}(\b C_B,\b Z_B) \sim \frac{1}{2} \frac{|f^{\prime}(0-)|}{f(0)}
 \frac{\BE|R - \BE R|}{s_n},\ {\rm\ if\ }f^\prime(0-) \neq 0,
 \en
and
 \eq\label{approx4gen}
 d_{TV}(\b C_B,\b Z_B) = o\left(\frac{1}{s_n}\right),\ {\rm\ if\ }
f^\prime(0-) = 0.
 \en
For partitions of an integer and for partitions of a set, a good choice
for $s_n$ is the standard deviation $\sigma_n$
 with asymptotics given by (\ref{vartn-multipart}) and
(\ref{Moser2}), and $X$ is normally distributed, so that
(\ref{approx4gen}) should apply.

 Observe that for two fixed sets $B, B^{\prime}$ the approximation  in
 (\ref{approx3}) or (\ref{approx3gen}) has as a corollary the
statement  that if $f^\prime(0-) \neq 0$ then
as $n \to \infty$,
 $$
 \frac{d_{TV}(\b C_B,\b Z_B)}{d_{TV}(\b C_{B^{\prime}},\b Z_{B^{\prime}})} \to
 \frac{\BE |R_B - \BE R_B|}{ \BE |R_{B^{\prime}} - \BE R_{B^{\prime}}|}.
 $$
 
 By the Cauchy-Schwarz inequality,
 $\BE |R_B - \BE R_B| \leq \sigma_B$, so another rigorous version
 of the heuristic in (\ref{approx2}) would be the statement that as $n \to
 \infty$, $d_{TV}(\b C_B, \b Z_B) = O( \sigma_B/ \sigma_n)$ uniformly in
 $B$; that is
 \eq\label{approx5}
 \lim_{n \to \infty} \sup_{B \subset [n]} \left( d_{TV}(\b C_B, \b Z_B)
 \frac{\sigma_n}{\sigma_B}\right) < \infty.
 \en
 Note that (\ref{approx5}) is not embarrassed by the largest possible
 $B$, namely $B = [n]$, since $d_{TV}(\cdot,\cdot) \leq 1$.
 
 
\subsection{Examples with a limit process: the logarithmic class}
\label{logsect}
 
The previous section suggests that the limit law of $T_n / n$ plays a key
role in analyzing the accuracy  of the approximation of certain
combinatorial structures by independent processes. 
The logarithmic class consists of those assemblies which satisfy
(\ref{assembly-hyp}), and those multisets and selections which satisfy
(\ref{multiset-hyp}).  All of these, with the appropriate constant
choice of $x$, satisfy 
\eq\label{Zcond}
i \B EZ_i \rightarrow \kappa,\ i \BP(Z_i=1)
\rightarrow \kappa {\rm \ for\ some\ } \kappa >0.
\en
Lemma \ref{tnn-mul} below shows that, for $Z_i$ satisfying
(\ref{Zcond}), and $T_n=Z_1+2Z_2+\cdots+nZ_n$, the limit distribution of
$T_n /n$ depends only on the parameter $\kappa$.

Let $d_W$ be the $L_1$ Wasserstein distance
between distributions, which  can be defined, in the same
spirit as (\ref{tvcouple}), by 
$$
d_W(X,Y) =
\min_{couplings} \BE |X-Y|.
$$
For $\BZ^+$-valued random variables,
$d_W$ is easily computed via 
$$
d_W(X,Y)=\sum_{i \geq 1} |\BP(X\geq i)-\BP(Y\geq i)|,
$$
and when $X$ is stochastically larger than $Y$, so that the absolute
values above do nothing, this further simplifies to $d_W(X,Y)=\BE X-\BE
Y$. Note that for integer-valued random values, $d_W \geq d_{TV}$.  

Let  $\tilde Z_i$ be
 Bernoulli with parameter $\kappa/i \wedge 1$, and let 
$Z_i^*$ be Poisson with mean $\kappa/i$. It is easy to check that the
 condition (\ref{Zcond}) is equivalent to 
$d_W(Z_i,\tilde Z_i)=o(1/i)$. Since $d_W(\tilde Z_i,Z_i^*)=o(1/i)$, the triangle
inequality implies that
the condition (\ref{Zcond}) is also equivalent to 
$d_W(Z_i,Z_i^*)=o(1/i)$. 

For the class of assemblies that satisfy the condition
(\ref{assembly-hyp}), we use $x=y^{-1}$ and $\BE Z_i=m_i x^i/i!$, so
that $\BE Z_i \sim \kappa/i$.  Lemma \ref{tnn-ass} applies directly;
for Poisson random variables (\ref{Zcond}) is equivalent to $\B E Z_i
\sim \kappa/i$, so Lemma \ref{tnn-mul} also applies. 
For multisets  and selections 
satisfying the hypothesis (\ref{multiset-hyp}), it is easy to show
that (\ref{Zcond}) holds.

\begin{lemma}\label{tnn-ass}
If $Z_j$ are
independent Poisson random variables with $\BE Z_j = \lambda_j \sim
\kappa /j$ for some constant $\kappa >0$, and
$T_n = \sum_{j=1}^n j Z_j$, then 
\eq\label{limitlaw}
n^{-1} T_n \Rightarrow X_\kappa,\ n \to \infty
\en
and $X_\kappa$ has Laplace transform 
\eq\label{limitlt}
\psi(s) \equiv \BE e^{- s X_\kappa} = \exp\left( - \kappa \int_0^1 
(1 - e^{-s x}) \frac{dx}{x} \right).
\en
\end{lemma}

\proof By direct calculation,
\begin{eqnarray*}
\log \BE e^{- s T_n / n} & = & - \sum_{j = 1}^n \lambda_j (1 -
e^{-j s / n}) \\
&&\\
& = & - \sum_{j = 1}^n \frac{\kappa}{j} (1 - e^{-j s / n}) +
\sum_{j = 1}^n \left( \frac{\kappa}{j} - \lambda_j \right) (1 -
e^{-j s / n})
\end{eqnarray*}
Clearly, the first term on the right converges to $-\kappa \int_0^1(1-
e^{-sx})\frac{dx}{x}$. That
the second term is $o(1)$ follows by observing that $\lambda_j -
\kappa / j = o(j^{-1})$, and comparing to the first sum.
\hfill \qed

\begin{lemma}\label{tnn-mul}
For $i=1,2,\ldots$, let $Z_i$ be positive integer-valued random
variables satisfying the conditions in (\ref{Zcond}).
If $T_n = \sum_{j=1}^n j Z_j$, then 
\eq\label{limitlawmul}
n^{-1} T_n \Rightarrow X_\kappa,\ n \to \infty
\en
and $X_\kappa$ has Laplace transform given in (\ref{limitlt}).
\end{lemma}

\proof Construct independent Bernoulli random variables $\tilde Z_i =
Z_i \wedge 1$. Clearly $\tilde Z_i \leq Z_i$ and
$\BP(Z_i = 1) \leq \BE \tilde Z_i \leq \BE Z_i$. It follows that 
$i \BE \tilde Z_i \to \kappa$. Therefore
$$
i|\BE Z_i - \BE \tilde Z_i| = i(\BE Z_i - \BE \tilde Z_i) \to 0.
$$
Hence if $\tilde T_n = \tilde Z_1 + \cdots + n \tilde Z_n$,
\eq\label{useme}
\BE \left|\frac{T_n}{n} - \frac{\tilde T_n}{n} \right| \to 0.
\en
It remains to show that $n^{-1} \tilde T_n \todist X_{\kappa}$.

For $i=1,2,\ldots$, let $Z_i^*$ be independent Poisson random variables
satisfying $p_i \equiv \BE Z_i^* = \BE \tilde Z_i \sim \kappa/i.$ 
We may construct
$Z_i^*$ in such a way that for each $i$
$$
\BE |\tilde Z_i - Z_i^*| = d_W(\tilde Z_i, Z_i^*),
$$
where $d_W$ denotes Wasserstein $L_1$ distance. But if $X$ is Bernoulli
with parameter $p$ and $Y$ is Poisson with parameter $p$, then a simple
calculation shows that $d_W(X,Y) = 2(p - 1 + e^{-p}) \leq p^2$.
Hence 
$$
n^{-1} \BE|\tilde T_n - T_n^*| \leq n^{-1} \sum_{i=1}^n i p_i^2 \to 0.
$$
It follows that $n^{-1} \tilde T_n$ has the same limit law as $n^{-1}
T_n^*$, which is that of $X_{\kappa}$ by Lemma
\ref{tnn-ass}. \hfill\qed

The random variable $X_\kappa$ has appeared in several guises before, not
least as part of the description of the density of points in a 
Poisson--Dirichlet process. See Watterson (1976), Vershik and Shmidt
(1977), Ignatov (1982), Griffiths (1988) and Ethier and Kurtz (1986) and 
the references contained therein. 
For our purposes, it is enough to record that the density $g(\cdot)$ of
$X_\kappa$ is known explicitly on the interval $[0,1]$:
\eq\label{ignatov}
g(z) = \frac{e^{- \gamma \kappa}}{\Gamma(\kappa)} z^{\kappa - 1},\ 0
\leq z \leq 1,
\en
where $\gamma$ is Euler's constant.
From (\ref{ignatov}) follows the fact that
\eq\label{ignatov1}
\frac{g^\prime(1-)}{g(1)} = \kappa - 1.
\en
We may now combine the previous results with (\ref{approx2}) and
(\ref{mean-varass}) to rephrase the asymptotic
behavior of $d_{TV}(\bC_B,\bZ_B)$ in (\ref{approx3}) and
(\ref{approx4}) as follows. For any assembly satisfying
(\ref{assembly-hyp}), or for any multiset or selection satisfying
(\ref{multiset-hyp}), we should have the following decay rates, 
for any fixed $B$, as $n \rightarrow \infty$.

 {\it In the case $\kappa \neq 1$}
 \eq\label{approx3b}
 d_{TV}(\b C_B,\b Z_B) \sim \frac{1}{2} |\kappa - 1|
 \frac{\BE|R - \BE R|}{n},
 \en
 
 {\it In the case $\kappa = 1$}
 \eq\label{approx4b}
 d_{TV}(\b C_B,\b Z_B) = o\left(\frac{1}{n}\right).
 \en

 For a class of examples known as the Ewens
sampling formula, described in Section \ref{esfsect},
and for $B$ of the form 
$B = \{1,2,\ldots,b\}$, (\ref{approx3b}) is
proved in Arratia, Stark and Tavar\'e (1994). The analogous result for 
random mappings, in
which   $\kappa = 1/2$, and other assemblies that can be approximated by the
Ewens sampling formula, may also be found there.
For the corresponding results for multisets and selections, see Stark
(1994b).

The statement (\ref{approx4b}) has been established for random
permutations by Arratia and Tavar\'e (1992), where it is shown inter
alia that for $B = \{1,2,\ldots,b\}$, $d_{TV}(\bC_B,\bZ_B) \leq F(n/b)$,
where $\log F(x) \sim -x \log x$ as $x \to \infty$. For the case of
random polynomials over a finite field, Arratia, Barbour and Tavar\'e
(1993) established that $d_{TV}(\bC_B,\bZ_B) = O(b \exp(-c n / b))$,
where $c = \frac{1}{2}\log(4/3)$.

Among the class of assemblies in the logarithmic class, weak convergence 
(in $\BR^\infty$) of the component counting process to
the appropriate Poisson process has been established for random
permutations by Goncharov (1944), for random mappings by Kolchin (1976),
and for the Ewens sampling formula by Arratia, Barbour and Tavar\'e
(1992). For multisets in the logarithmic class, this has been
established for random polynomials by Diaconis, McGrath and Pitman
(1994) and Arratia, Barbour and Tavar\'e (1993), and for random mapping
patterns by Mutafciev (1988).

  \section{Non-uniqueness in the choice of the parameter $x$}
 \label{sect5}

 An appropriate choice of $x = x(n)$ for good approximation is not
 unique. 
An obvious candidate is that $x$ which maximizes $\BP(T_n = n)$, which
is also that $x$ for which $\BE T_n = n$. This can be seen by
differentiating $\log \BP(T_n = n)$ in formulas (24) - (26) and
comparing to $\BE T_n$ from formulas (11) - (13); at the general level
this is the observation that $\BP(T=t)$ in (19) is maximized by that $x$
for which $\BE T = t$. Nevertheless, the obvious candidate is not always
the best one.
We discuss here two qualitatively different examples: the
logarithmic class, and partitions of a set.

\subsection{The Ewens sampling formula}\label{esfsect}

The central object in the logarithmic class is the Ewens sampling
formula (ESF). This is the family of distributions with parameter
$\kappa > 0$ given by
(\ref{equaldist}), where the $Z_i$ are independent Poisson random
variables with $\BE Z_i = \kappa/i$, or more generally, with
\eq\label{esfmeans}
\lambda_i \equiv \BE Z_i = \frac{\kappa x^i}{i},
\en
the conditional distribution being unaffected by the choice of $x >
0$.
For $\kappa = 1$, the ESF is the distribution of cycle counts for a
uniformly chosen random permutation. For $\kappa \ne 1$, the ESF 
can be viewed as the nonuniform measure
on permutations with sampling bias proportional to $\kappa^{{\rm \#\
cycles}}$; see Section \ref{sect8} for details. 
The ESF arose first in the context
of population genetics (Ewens, 1972), and is given explicitly by
\eq\label{esfdef}
\BP(C_1(n)=a_1,\ldots,C_n(n) = a_n) =
\bone(\sum_{l=1}^n l a_l = n) \ \frac{n!}{\kappa_{(n)}} \ \prod_{i=1}^n
\left( \frac{\kappa}{i} \right)^{a_i} \frac{1}{a_i!} .
\en
The ESF corresponds to (\ref{assembly-hyp}) with $y = 1$ and
the asymptotic relation in $i$ replaced by equality. It is useful in
describing all assemblies, multisets and selections in the logarithmic
class; see Arratia, Barbour and Tavar\'e (1994) for further details.

For irrational $\kappa$ the ESF cannot  be realized as
a uniform measure on a class of combinatorial objects.  For rational
$\kappa =r/s$ with integers $r>0, s>0$, there are at least two
possibilities.  First, comparing  
(\ref{equaldist}) with $\B EZ_i = \kappa/i$, and (\ref{Zassembly})
with $\B EZ_i= m_i x^i /i!,$ for any choice $x>0$, we take $x=1/s$ to 
see that the ESF
is the uniform measure on the assembly with 
$m_i =r(i-1)!s^{i-1}$.  One interpretation of this is permutations on
integers, enriched by coloring each cycle with one of $r$ possible
colors, and coloring each element of each cycle, except the smallest,
with one of $s$ colors.  For a second construction, we use a device from
Stark (1994a).  Consider permutations of $ns$ objects, 
in which all
cycle lengths must be multiples of $s$.  Formally, this is the assembly
on $ns$ objects, with $m_i=(i-1)!\bone (s|i)$, so that
$(C_1,C_2,\ldots,C_{ns}) \indist
(Z_1,Z_2,\ldots,Z_{ns}|Z_1+2Z_2+\cdots+nsZ_{ns}=ns)$, where $Z_i$ is
Poisson with $\B EZ_i= \bone(s|i) \ /i$.  Since those $C_i$ and
$Z_i$ for which $s$ does not divide $i$ are identically zero, we
consider $C_i^*\equiv C_{is}, \ Z_i^* \equiv Z_{is}$, and $T_n^* \equiv
Z_1^*+2Z_2^*+\cdots+nZ_n^* = \frac{1}{s}(Z_1+2Z_2+\cdots+nsZ_{ns})$. 
We have $(C_1^*,\ldots,C_n^*) \indist (Z_1^*,\ldots,Z_n^* | T_n^* = n)$, and the
$Z_i^*$ are independent Poisson with $\B E Z_i^*=1/(si)$.  Thus
the distribution of $(C_1^*(n),\ldots,C_n^*(n))$ is the ESF with
$\kappa=1/s$.  To change this to $\kappa=r/s$, we need only
color each cycle with one of $r$ possible colors, so that $m_i = r (i-
1)! \bone(s|i), \B EZ_i = r \bone (s/i)\ /i$, and $\B EZ_i^* =
r/(si)$.  To summarize our second construction of the ESF with
$\kappa =r/s$, let $C_i^*(n)$ be the number of cycles of length $si$ in
a random permutation of $ns$ objects, requiring that all cycle lengths
be multiples of $s$, and assigning one of $r$ possible colors to each
cycle. 

For comparing the ESF to the unconditioned, independent process
$(Z_1,\ldots,$ $Z_n)$ it is interesting to consider the role of varying $x$.
 The choice $x = 1$ in (\ref{esfmeans}), so that $\BE Z_i = \kappa /i$,
yields $\BE T_n = \kappa n$, and $\sigma_n \sim n \sqrt{\kappa/2}$.
 In the case $\kappa \ne 1$ the discrepancy between $\BE T_n$ and the
 goal $n$ is a bounded multiple of $\sigma_n$. This is close enough
 for good approximation, in the sense that
 $(C_1(n),\ldots,C_n(n),0,\ldots)$ $\Rightarrow (Z_1,Z_2,\ldots)$.
 This, together with a $O(b/n)$ bound on
 $d_{TV}(C_1(n),\ldots,C_b(n)),(Z_1,\ldots,Z_b))$ that is  uniform in $1 \leq b
 \leq n$, is proved in \TVthree \ by exploiting a coupling
 based on Feller (1945). This coupling
 provided even stronger information whose utility is discussed in \TVtwo.
 Barbour (1992) showed that the $O(b/n)$ bound above cannot be replaced
 by $o(b/n)$ for $x=1$, $\kappa \ne 1$.

 For the case of independent $Z_i^{\prime}$ which are Poisson with means
varying with $n$ given by
 $$
 \BE Z_i^{\prime} = \BE C_i(n) = \frac{\kappa}{i} \frac{n (n-
 1)\cdots(n-i+1)}{(\kappa+n-i)\cdots(\kappa+n-1)},
 $$
 Barbour (1992) showed that
 $d_{TV}(C_1(n),\ldots,C_b(n)),(Z_1^{\prime},
 \ldots,Z_b^{\prime})) =
 O((b/n)^2)$, uniformly in $1 \leq b \leq n$. Observe that with this
 choice of Poisson parameters, $\BE T^\prime_n \sim \kappa n$ but it is 
{\it not} 
the case 
that $(C_1(n),\ldots,C_n(n)) \indist (Z_1^\prime,\ldots,Z_n^\prime | 
T^\prime_n = n)$.

If we are willing to use coordinates $Z_i \equiv Z_i(n)$ whose means
vary with $n$, we can still have the conditional relation
(\ref{equaldist}) by using $x = x(n)$ in (\ref{esfmeans}). An
appealing family of choices is given by $x = \exp(- c/n)$, since this
yields for $c \neq 0$
\eq\label{newmean}
\BE T_n = \sum_{i=1}^n i \lambda_i = \sum_{i=1}^n i \frac{\kappa}{i}
e^{-i c / n} \sim n \frac{\kappa(1 - e^{-c})}{c}.
\en
 By choosing $c \equiv
 c(\kappa)$ as the solution of $\kappa = c/(1-e^{-c})$, we can make $\BE
 T_n \sim n$, and this should provide a closer approximation than the
choice $c=0, x=1$. 
 However an even better choice of $c$ is available. We explore this in
the next section.
 
 \subsection{More accurate approximations to the logarithmic class}
 \label{sect5.2}

For assemblies, multisets, and selections in the logarithmic class
discussed in Section \ref{logsect}, as
well as for the ESF, the choice  of $x$ proportional to $\exp(-c/n)$ is
interesting. In this situation, the limit law of $T_n/n$ depends only on
the parameters $\kappa$ and $c$. Properties of this limit law lead to an
optimal choice for $c$.

The following lemma applies to assemblies that satisfy the condition
(\ref{assembly-hyp}), and to the ESF by taking $m_i = \kappa (i-1)!,
y=1$, the $m_i$ not necessarily being integers.

\begin{lemma}\label{tnn-assnew}
Assume that $m_i \geq 0$ satisfies 
$m_i/i! \sim \kappa y^i / i$ for constants $y \geq 1, \kappa > 0$, 
and set $x = e^{- c / n} y^{-1}$ for constant $c \in \BR$. If $Z_j
\equiv Z_j(n)$ are
independent Poisson random variables with $\BE Z_j = m_j x^j / j!$, and
$T_n = \sum_{j=1}^n j Z_j$, then 
\eq\label{limitlaw1}
n^{-1} T_n \Rightarrow X_{\kappa,c},\ n \to \infty
\en
and $X_{\kappa,c}$ has Laplace transform 
\eq\label{limitlt1}
\psi_c(s) \equiv \BE e^{- s X_{\kappa,c}} = \exp\left( - \kappa \int_0^1 
(1 - e^{-s x}) \frac{e^{-c x}}{x} dx \right).
\en
\end{lemma}

\proof As in Lemma \ref{tnn-ass}, calculate the limit of the 
log Laplace transform.\hfill \qed

Next we prove that the same limit law holds 
for multisets or selections satisfying the hypothesis (\ref{multiset-hyp}). 
\begin{lemma}\label{tnn-mulnew}
Assume that the multiset (or selection) satisfies (\ref{multiset-hyp}):
$m_i \sim \kappa y^i / i$ for constants $y \geq 1, \kappa >0$,
 and set $x = e^{ - c / n} y^{-1}$. If $Z_j
\equiv Z_j(n)$ are
independent negative binomial random variables with parameters $m_j$ and
$ x^j$ (respectively, binomial with parameters $m_j$ and $x^j/(1+x^j))$ 
and $T_n = \sum_{j=1}^n j Z_j$, then 
\eq\label{limitlawmul2}
n^{-1} T_n \Rightarrow X_{\kappa,c},\ n \to \infty
\en
and $X_{\kappa,c}$ has Laplace transform given in (\ref{limitlt1}).
\end{lemma}

\no{\bf Remark: } For the case of multisets, we assume that $x < 1$.

\proof Observe first that in either case, if $b = o(n)$, then 
$n^{-1} \BE T_{0b} \rightarrow 0$, so that $n^{-1} T_{0b} \inprob 0$ as
$n \to \infty$. Let $\tilde Z_j$ be independent Poisson random variables
with $\BE \tilde Z_j = m_j x^j$, and write 
$\tilde T_n = \sum_{j=1}^n j
\tilde Z_j$, $\tilde T_{bn} = \sum_{j=b+1}^n j \tilde Z_j$. We show
that for $b = o(n)$, $T_{bn}/n$ and $\tilde T_{bn} / n$ have the same limit 
law, which completes the proof since by Lemma \ref{tnn-assnew},
$\tilde T_{bn} / n \Rightarrow X_{\kappa,c}$.
We will use the notation NB, Po, and Geom to denote the negative
binomial, Poisson and geometric distributions with the indicated
parameters.

For the multiset case, notice that
\begin{eqnarray*}
d_{TV}(T_{bn},\tilde T_{bn}) & \leq & d_{TV}((Z_{b+1},\ldots,Z_n),
(\tilde Z_{b+1},\ldots,\tilde Z_n)) \\
& \leq & \sum_{b+1}^n d_{TV}(Z_j,\tilde Z_j).
\end{eqnarray*}
To estimate each summand, we have
\begin{eqnarray}
d_{TV}(Z_j, \tilde Z_j) & = & d_{TV}({\rm NB}(m_j, x^j), {\rm Po}(m_j
x^j)) \nonumber \\
& \leq & m_j d_{TV}( {\rm Geom}(x^j), {\rm Po}(x^j)) \nonumber\\
& \leq & 2 m_j x^{2j}. \label{qandd1}
\end{eqnarray}
The bound in (\ref{qandd1}) follows from the fact that
$d_{TV}({\rm Geom}(p), {\rm Be}(p)) = p^2$ and  $d_{TV}({\rm Be}(p), {\rm
Po}(p)) = p(1-e^{-p}) \leq p^2$, so that
$d_{TV}({\rm Geom}(p),{\rm Po}(p)) \leq 
d_{TV}({\rm Geom}(p), {\rm Be}(p))$ $+\ d_{TV}({\rm Be}(p), {\rm Po}(p)) 
\leq 2p^2$, a result we apply with $p = x^j$. Hence
$$
d_{TV}(T_{bn}, \tilde T_{bn}) \leq 2 \sum_{j=b+1}^n (m_j x^j) x^j =
O(y^{-b}/b).
$$
Choosing $b \to \infty, b = o(n)$ completes the proof for multisets.

For the selection case, (\ref{qandd1}) may be replaced by
$$
d_{TV}(Z_j, \tilde Z_j) \leq m_j d_{TV}({\rm Be}(x^j/(1 + x^j)), {\rm
Po}(x^j)) \leq 2 m_j x^{2j}.
$$
The last estimate following from the observation that $d_{TV}({\rm
Be}(p/(1 + p)), {\rm Be}(p)) = p^2/(1+p)$, so that
$d_{TV}({\rm Be}(p/(1 + p)),{\rm Po}(p)) \leq d_{TV}({\rm
Be}(p/(1 + p)), {\rm Be}(p))$ $ + d_{TV}({\rm Be}(p), {\rm Po}(p)) 
\leq 2p^2$, which we apply with $p = x^j$. This completes the proof. 
\hfill \qed

The random variable $X_{\kappa}$ of Section \ref{logsect} is the special case
$c = 0$ of $X_{\kappa,c}$. Further, for $c \neq 0$,
$$
\BE X_{\kappa,c} = \kappa \frac{1 - e^{-c}}{c}
$$
and
$$
{\rm Var} X_{\kappa,c} = \kappa \frac{1 - (1 + c)e^{-c}}{c^2}.
$$
The density $g_c$ of $X_{\kappa,c}$ may be found from the density $g$ of
$X_{\kappa}$ by observing that the log Laplace transforms, given by
(\ref{limitlt}) and (\ref{limitlt1}), are related by
$$
\psi_c(s) = \frac{\psi(c+s)}{\psi(c)}
$$
so that
$$
g_c(z) = e^{-c z} g(z) / \psi(c), \ z \geq 0.
$$
In particular, from (\ref{ignatov}),
\eq\label{gcdensity}
g_c(z) = \frac{e^{- \gamma \kappa} e^{- c z} z^{\kappa -
1}}{\Gamma(\kappa) \psi(c)},\ 0 \leq z \leq 1.
\en

From (\ref{gcdensity}) the
value of $c$ that maximizes the density
$g_c(z)$ for fixed $z \in [0,1]$ is the $c$ that maximizes
$-c z - \log \psi(c)$, just as suggested by large deviation theory. 
This $c$ is the solution of the equation
$$
c z = \kappa ( 1 - e^{- c z}).
$$
Using $z = 1$, we see from the heuristic (\ref{local1}) that choosing
$c$ to be the
 solution of $c= \kappa (1-e^{-c})$ asymptotically maximizes $\BP(T_n = n)$; 
and from (\ref{newmean}), this  also makes $\BE T_n \sim n$.

However, the heuristic in (\ref{approx3}) and 
(\ref{approx4}) suggests that better
approximation should follow from choosing $c$ so that $g_c^{\prime}(1-) =
0$. From (\ref{gcdensity}) and (\ref{ignatov1}), we get 
\eq\label{ignatov3}
c = \frac{g^{\prime}(1 -)}{g(1)} = \kappa - 1.
\en
For this choice of $c$ we have $g_c^{\prime}(1-) = 0$, and 
\eq\label{ignatov4}
\frac{g_c^{\prime\prime}(1 -)}{g_c(1)} = 1 - \kappa.
\en
A second order approximation in the spirit of Section \ref{sect4} then leads 
us to the following heuristic: for any fixed $B$,

 {\it In the case $\kappa \neq 1$}
 \eq\label{approx3c}
 d_{TV}(\b C_B,\b Z_B) \bothsides \frac{\sigma_B^2}{n^2},
 \en
 
 {\it In the case $\kappa = 1$}
 \eq\label{approx4c}
 d_{TV}(\b C_B,\b Z_B) = o\left(\frac{1}{n^2}\right).
 \en

For the case $B = [b] \equiv \{1,2,\ldots,b\}$, extensive numerical
computations using the recurrence methods described in Section \ref{sect9}
support these conjectures for several of the combinatorial examples
discussed earlier. In these cases, the bound in (\ref{approx3c}) is of
order $(b/n)^2$. Finding the asymptotic form of this rate seems to be a
much harder problem, since it seems to depend heavily on the value of
$\kappa$.

\subsection{Further examples}\label{sect5.3}

The class of partitions of a set provides another example to show that 
the choice of
 $x$ for good approximation is partly a matter of taste. In this example,
 $m_i \equiv 1$, so that
 $$
 \BE T_n = \sum_{i=1}^n \frac{i m_i x^i}{i!} = x \sum_{i=0}^{n-1}
 \frac{x^i}{i!}.
 $$
 One choice of $x$ would be the exact solution $x^*$ of the equation $\BE
 T_n = n$, but this choice is poor since the definition of $x^*$ 
is complicated. 
 A second choice which is more usable is to take $x = x^\prime$, the 
 solution of the equation $x e^x = n$. This is based on the
 observation that $\BE T_n \sim x e^x$, provided $x = o(n)$.
 The solution  $x^\prime$ has the form
 (cf. de Bruijn, 1981, p. 26)
 $$
 x^\prime = \log n - \log \log n + \frac{\log \log n}{\log n} + \frac{1}{2}
 \left( \frac{\log \log n}{\log n} \right)^2 + O\left( \frac{\log \log
 n}{\log^2 n} \right).
 $$

 For set partitions,
 with either $x^*$ or $x^\prime$ in the role of $x$, we have $\sigma_n^2
  \sim x^2 e^x \sim n \log n$, and we can check that $| n - \BE T_n| =
 O(\sqrt{n \log n})$ is satisfied using $x=x^\prime$. This corresponds
to checking the condition in (\ref{metathm1}). Comparing the condition
$\BE T_n \sim n$ with the condition that $n - \BE T_n = O(\sigma_n)$
required by (\ref{metathm1}), we see that in the logarithmic class the
former is too restrictive while for set partitions it is not
restrictive enough.

 \section{Refining the combinatorial and independent
processes}\label{sect6}

\subsection{Refining and conditioning}

Although the refinements considered in this section are complicated in
notation, the ingredients -- including geometric and Bernoulli random
variables and the counting formulas (\ref{Rassembly}) -
(\ref{Rselection}) -- are
simpler than their unrefined counterparts.

The dependent random variables $C_i \equiv C_i(n)$, which
 count the number of
components of weight $i$ in a randomly selected object of total
weight $n$, may be refined as
\[
  C_i = \sum_{j = 1}^{m_i} D_{ij}.
\]
Here we suppose that the $m_i$ possible structures
of weight $i$ have been labelled $1,2,\ldots,m_i$, and $D_{ij} \equiv
D_{ij}(n)$ counts the number of occurrences of the $j^{th}$ object of
weight $i$.  The independent random variable $Z_i$ can also be refined, as
\[
  Z_i = \sum_{j=i }^{ m_i} Y_{ij},
\]
where the $Y_{ij}$ are mutually independent, and for each $i$,
$Y_{i1},Y_{i2},\ldots,Y_{im_i}$ are identically distributed.  For
assemblies, multisets, and selections respectively, the
distribution of $Y_{ij}$ is Poisson $( x^i/i!)$ for $x>0$,
 geometric($x^i$) for $0<x<1$, or
Bernoulli($x^i/(1+x^i))$ for $x>0$. If the choice of parameter $x$ is
taken as a function of $n$, then one can view $Y_{ij}$ as $Y_{ij}(n)$.
For assemblies, with $x > 0$,
\eq\label{80a}
\BP(Y_{ij}=k) = \exp(-x^i/i!) \frac{(x^i/i!)^k}{k!},\ k=0,1,\ldots.
\en
For multisets, with $0 < x < 1$, 
\eq\label{Zrefine}
\BP(Y_{ij}=k) = (1-x^i)x^{ik},\ k=0,1,\ldots,
\en
whereas for selections, with $x > 0$, we have
\eq\label{80c}
\BP(Y_{ij}=k) = 
 \frac{1}{1+x^i} \bone(k=0) \ + \ \frac{x^i}{1+x^i} \bone(k=1).
\en

For the full refined processes corresponding to a random object of size
$n$ we denote the combinatorial process by
\[
\b D(n) \equiv (D_{ij}(n),\ 1 \leq i \leq n, 1 \leq j \leq m_i),
\]
and the independent process by
\[
\b Y(n) \equiv (Y_{ij},\ 1 \leq i \leq n, 1 \leq j \leq m_i).
\]
The weighted sum $T_n = \sum_1^n i Z_i$ is of course a
weighted sum of the refined independent $Y$'s, since
\[
  T_n = \sum_{i=1}^n \sum_{j=1}^{m_i} i Y_{ij}.
\]
\begin{theorem}\label{refinedthm}
For assemblies, multisets, and selections, if $\BP(T_n=n)>0$, then 
the refined combinatorial
process, for a uniformly chosen object of weight $n$, is equal in
distribution to the independent process $\b Y(n)$, conditioned on the
event $\{ T_n=n \}$, that is
\[ \b D(n) \indist (\b Y(n)|T_n=n).  \]
\end{theorem}
\proof
Just as (\ref{equaldist}) is a special case of Theorem
\ref{genequaldist} with $t=n$, so is this.  Imagine first the special case of
(\ref{equaldist}) with each $m_i \equiv 1$, and then replicate $m_i$--
fold the
index $i$ and its corresponding function $g_i$ and normalizing constant
$c_i$.  The case $m_i=0$ for some $i$ is allowed. We have index set
\eq
\label{refineI}
  I = \{ \alpha =(i,j): 1 \leq i \leq n, 1 \leq j  \leq m_i \}
\en
and weight function $w$ given by $ w(\alpha) = i$ for  $\alpha = (i,j) \in I$.

The reader should be convinced by now, but for the record, here are the
details.
For $\b{b} \equiv (b(\alpha))_{\alpha \in I} \in \B{Z}_+^I$,
write $\b b \cdot \b w \equiv \sum_I w(\alpha)b(\alpha)$. Consider the number 
$R(n,\b{b})$ of
objects of total weight $\b b \cdot \b w =n$, having $b_\alpha \equiv b(\alpha)$
 components of type $\alpha$, for $\alpha \in I$.
  For assemblies, the refined generalization  of Cauchy's formula is that
\begin{eqnarray}
R(n, \b{b}) & \equiv &  | \{ \mbox{assemblies on [n]}: \b{D} = \b{b} \} |
\nonumber \\
& = &
 \bone(\b b \cdot \b w =n) 
\   n! \ \prod_{\alpha \in I} \frac{1}{(i!)^{b(\alpha) } \ b(\alpha) !},
\label{Rassembly}
\end{eqnarray}
where $i =w(\alpha)= $ the first coordinate of $\alpha$. 
For multisets, 
\begin{eqnarray}
R(n, \b{b}) & \equiv &  | \{ \mbox{multisets of weight n}: \b{D} = \b{b} \} |
 \\ \nonumber
& = & \bone(\b b \cdot \b w =n),
\label{Rmultiset}
\end{eqnarray}
while for selections, 
\begin{eqnarray}
R(n, \b{b}) & \equiv &  | \{ \mbox{selections of weight n}: \b{D} = \b{b} \} |
\\  \nonumber 
& = &
 \bone(\b b \cdot \b w =n) 
\    \ \prod_1^n {1  \choose b_\alpha }.
\label{Rselection}
\end{eqnarray}

These examples have the form
\eq
R(n, \b{b})  \equiv   | \{  \b{D} = \b{b} \} | = 
 \bone(\b b \cdot \b w =n) 
\  f(n)  \ \prod_{\alpha \in I} g_\alpha(b_\alpha),
\label{Rgeneral}
\en
with $f(n)=n!$ for assemblies and $f(n) \equiv 1$ for multisets and
selections. With $p(n)$ given by (\ref{def p(n)}), we have the refined
analysis of the total number of structures of weight $n$:
\eq
   p(n)= \sum_{\b b \in \B Z_+^I} R(n,\b b).
\en
Picking an object of weight $n$ uniformly defines the refined
combinatorial distribution
\eq
  \B P(\b D(n) = \b b) \equiv   \frac{R(n,\b{b})}{p(n)} =
 \bone(\b b \cdot \b w =n) 
\  \frac{f(n)}{p(n)}  \ \prod_I g_\alpha(b_\alpha).
\label{refinedcombdist}
\en
Observe that with multisets, $g_\alpha(k)=1$ for $ k \in \BZ_+$;
with selections $g_\alpha(k) = {1 \choose k} = \bone(k=0 \mbox{ or }1)$;
and with assemblies, if $\alpha =(i,j)$, then $g_\alpha(k) =(1/i!)^k/k!,$
for $k \in \BZ_+$. 
Now apply Theorem \ref{genequaldist} with $\bD_I$ in the role of $\bC_I$,
$Y_{ij} \equiv Y_\alpha$ in the role of $Z_\alpha$, and $t=n$. \hfill \qed

{\bf Remark}. It would be reasonable to consider (\ref{Rassembly})
through (\ref{Rselection}) as the basic counting formulas, with
(\ref{Nassembly}) through (\ref{Nselection}) as corollaries derived by
summing, and to consider the Poisson, geometric, and Bernoulli
distributions in (\ref{Zrefine}) as the basic distributions, with the
Poisson, negative binomial, and binomial distributions in
(\ref{Zassembly}) through (\ref{Zselection}) derived by convolution. 

\subsection{Total variation distance}

Since the refined combinatorial process $\b D(n)$ and the refined
independent process $\b Y(n)$ are related by conditioning on the value
of a weighted sum of the $Y$'s, Theorem \ref{tvthm} applies. For $K
\subset I$, where $I$ is given by (\ref{refineI}), write $\b D_K$ and
$\b Y_K$ for our refined processes, restricted to indices in $K$. Write 
\[
R_K' \equiv \sum_{\alpha \in K} w(\alpha) Y_\alpha, \ \ \ \ \
S_K' \equiv \sum_{\alpha \in I-K} w(\alpha) Y_\alpha,
\]
so that $T \equiv T_n = R_K' + S_K'$. 

\begin{theorem}\label{tvrefine}
\eq
d_{TV}(\b D_K,\b Y_K) \ = d_{TV}((R_K'|T=n),R_K').
\label{refinedtveq,general}
\en
\end{theorem}

\proof
This is a special case of Theorem \ref{tvthm}, with the
independent process $\b Y(n) \equiv \b Y_I$ playing the role of $\b Z_I$
and $\b D(n) \equiv \b D_I$ playing the role of $\b C_I$.  
Theorem  
\ref{refinedthm} 
is used to verify that the hypothesis 
(\ref{tvhyp}) is satisfied, in the form $\b D_I  \indist
(\b Y_I|T=n)$.
\hfill \qed

For the special case where $B \subset \{1, \ldots,n \}$ and $K=\{\alpha
= (i,j) \in I: i \in B \}$,  denote the restriction of the refined
combinatorial process, restricted to sizes in $B$,
 by $\b D_{B^*} \equiv \b D_K$,
so that
\[
\b D_{B^*} \equiv (D_{ij},\ i \in B, 1 \leq j \leq m_i),
\]
and similarly define $\b Y_{B^*}$.  In this special case, $R_K'
= R_B \equiv \sum_{i \in B} i Z_i$ is the weighted sum, restricted to B,
for the unrefined process, so (\ref{refinedtveq,general}) reduces to
\eq
d_{TV}(\b D_{B^*},\b Y_{B^*}) \ = d_{TV}((R_B|T=n),R_B).
\label{refinedtveq}
\en
Furthermore, by Theorem \ref{tvthm} applied to the unrefined case, with
$I=\{1,\ldots,n \}$ and $w(i)=i$ , we see that  $ d_{TV}((R_B|T=n),R_B)$ 
is equal to
 $d_{TV}(\b C_B,\b Z_B)$.

We have here a most striking example of the situation analyzed in
Theorem \ref{tv=thm}, where taking functionals doesn't change a total
variation distance.  Namely, there is a functional $g: \B Z_+^I
\rightarrow \B Z_+^n$, which ``unrefines'', and the functional $h: \B
Z_+^B \rightarrow \B Z_+$ discussed in our second proof of Theorem
\ref{tvthm}, such that 
\[
g(\b D_{B^*})=\b C_B, \ g(\b Y_{B^*})=\b Z_B, \ \ \ \
h(\b C_B) \indist (R_B|T=n) , \mbox{ and }h(\b Z_B)=R_B,
\]
so that, a priori via (\ref{tvinequality}),
\eq
d_{TV}(\b D_{B^*},\b Y_{B^*}) \geq d_{TV}(\b C_B,\b Z_B) \geq
d_{TV}((R_B|T=n),R_B).
\label{surprise=}
\en
Perhaps the result in (\ref{refinedtveq}), which shows that equality
holds throughout (\ref{surprise=}), is surprising.

 \section{Conditioning on events of moderate probability}\label{sect7}

We consider random combinatorial structures conditioned on some event.
 Given that we have a good
approximation by another process, this other
process, conditioned on the same event, may yield a good approximation to
the conditioned combinatorial structure.  The conditioning event must
have moderate probability, large relative to the original approximation
error. In contrast, if the conditioning event is very unlikely then the 
approximating process must also be changed, as discussed in Section
\ref{sect8} on large deviations. 

\subsection{Bounds for conditioned structures}

In this subsection, we consider bounds on total variation distance 
that are inherited from an existing approximation, after additional
conditioning is applied. 
\begin{theorem}\label{mild condition}
Let $A \subseteq B \subseteq
[n]$, and let  $h: \B Z_+^B \rightarrow \{0,1\}$ be measurable with
respect to coordinates in $A$. Let $\b Z_B$, and $\b C_B$ be arbitrary
processes with values in $\B Z_+^B$, and let $\b Z_A$ and $\b C_A $ denote
their respective restrictions to coordinates in $A$. Let
\[
\b C_B^* \indist (\b C_B | h(\b C_B) = 1),
\]
and
\[
\b Z_B^* \indist (\b Z_B | h(\b Z_B) = 1).
\]
Write
$p  = \B P(h(\b Z_B) = 1),$
$q  = \B P(h(\b C_B) = 1),$
$d_B =  d_{TV}(\b C_B, \b Z_B)$, 
$d_A  =  d_{TV}(\b C_A, \b Z_A)$, and
assume that $p>0$ and $q>0$. Then
\begin{eqnarray}\label{tvbnd0}
d_{TV}(\b C_B^*, \b Z_B^*) 
& \leq & \frac{1}{2} \left| 1 - \frac{q}{p} \right| + \frac{d_B}{p} \\
& & \nonumber \\
& \leq & \frac{1}{p} \left( \frac{d_A}{2} +
d_B \right) \label{tvbnd1} \\
& & \nonumber \\
\label{tvbnd2}
& \leq & \frac{3}{2} \frac{d_B}{p}.
\end{eqnarray}

\end{theorem}
\proof
 The second to last inequality follows from 
the relation $| p - q | \leq d_A,$
and is useful when this is the extent of our ability to estimate
$q$.
The last inequality follows simply from the fact that $d_A \leq d_B$. To
establish the first inequality, we have
\begin{eqnarray*}
d_{TV}(\b C_B^*, \b Z_B^*) & = & \frac{1}{2} \sum_{\smb a \in \B Z_+^B} |\B
P(\b C_B^* = \b{a}) - \B P(\b Z_B^* = \b{a}) |\\
&&\\
& = & \frac{1}{2} \sum_{\smb{a}: h(\smb{a}) = 1} \left| \frac{\B P(\b C_B
= \b{a})}{q} - \frac{\B P(\b Z_B = \b{a})}{p} \right| \\
&&\\
& = & \frac{1}{2} \sum_{\smb{a}: h(\smb{a}) = 1} \left| \B P(\b C_B
= \b{a}) \left( \frac{1}{q} - \frac{1}{p}\right) +
\frac{\B P(\b C_B
= \b{a}) - \B P(\b Z_B = \b{a})}{p} \right| \\
&&\\
& \leq & \frac{1}{2} \left| \frac{1}{q} - \frac{1}{p} \right| 
\sum_{\smb{a} : h(\smb{a}) = 1}\B P(\b C_B = \b{a}) \\
&&\\
& & \ \ \ + \frac{1}{2p} 
\sum_{\smb{a}: h(\smb{a}) = 1} \left| \B P(\b C_B
= \b{a}) - \B P(\b Z_B = \b{a}) \right| \\
&&\\
& = & \frac{1}{2} \left| \frac{1}{q} - \frac{1}{p} \right| q
+ \frac{1}{2p} 
\sum_{\smb{a}: h(\smb{a}) = 1} \left| \B P(\b C_B
= \b{a}) - \B P(\b Z_B = \b{a}) \right| \\
&&\\
& \leq & \frac{1}{2} \left| \frac{1}{q} - \frac{1}{p} \right| q + \frac{1}{2p} 
\sum_{\smb{a}} \left| \B P(\b C_B
= \b{a}) - \B P(\b Z_B = \b{a}) \right| \\
&&\\
& = & \frac{1}{2} \left| 1 - \frac{q}{p} \right| + \frac{d_B}{p}. 
\end{eqnarray*}
\hfill \qed

{\bf Remark.} 
While the theorem above uses the notation $C_B$ and $Z_B$ to suggest
applications where one process is obtained from an independent process by
conditioning, no such structure is required.  An arbitrary discrete space
$S$, together with an arbitrary functional $h:S \rightarrow \{0,1\}$, may
be encoded in terms of $S=\B Z_+^2$, with $A=\{1\}$ and $B=\{1,2\}$, so
that $h$ depends only on the first coordinate.  Thus Theorem \ref{mild
condition} applies to
discrete random objects in general.

\subsection{Examples}

\subsubsection{Random permutations}

In this case, the $Z_i$ are independent Poisson
distributed random variables, with $\lambda_i \equiv \B E Z_i = 1/i$.  
In \TVone \ it is
proved that for $1 \leq b \leq n$, the total variation distance $d_b(n)$
between $(C_1(n),\ldots,C_b(n))$ and $(Z_1,\ldots,Z_b)$ satisfies $d_b(n)
\leq F(n/b)$ where 
\begin{eqnarray}
  F(x) & \equiv &\sqrt{2 \pi m} \, \frac{2^{m-1}}{(m-1)!} \quad + 
\frac{1}{m!} + 3 \, \left( \d{\frac{x}{e}}\right)^{-x}, 
\quad \mbox{ with } m \equiv \lfloor x 
\rfloor \label{Fdef}\\
& & \nonumber\\
& \sim & \left( \frac{2e}{\lfloor x-1 \rfloor} \right)^{\lfloor x -1
\rfloor} \nonumber
\end{eqnarray}
as $x \rightarrow \infty$.
To get an approximation result for derangements, we use the functional $h$
having $h((a_1,\ldots,a_b))= \bone(a_1=0)$, with $A=\{1\}$ and
$B=\{1,2,\ldots,b\}$. This makes $\b C_B^*$ the process counting cycles of
size at most $b$ in a randomly chosen derangement, and $\b Z_B^*
=(Z_1^*,Z_2^*,\ldots,Z_b^*) \indist (0,Z_2,\ldots,Z_b)$. The total
variation distance $d_b^*(n)$ between $\b C_B^*$ and $\b Z_B^*$
satisfies $d_b^*(n) \leq (3/2)e \ F(n/b)$, simply by using
(\ref{tvbnd2}).  

Changing random permutations to random 
derangements is a special case of conditioning on some fixed conditions of
the form $C_i(n) =c_i, i \in A$, for given constants $c_i$, with $A \subseteq
B \subseteq \{1,2,\ldots,b\}$. In this situation, all the $Z_i^* $ are 
mutually independent, $Z_i^* \equiv c_i$ for $i \in A$,
and  for $i \notin A$,\ \ $Z_i^* \indist Z_i$ is Poisson with mean $1/i$.
   Here, Theorem \ref{mild condition}
yields the bound $d_b^*(n) \leq 3/(2p) F(n/b)$, where $p = \B P(Z_i=c_i
 \ \forall i \in A)$.  Theorem 3 in \TVone  \ gives a different upper bound, namely
$d_b^*(n) \leq F((n-s)/b) +2be ((n-s)/(be))^{-(n-s)/b}$, where $s= \sum_{i
\in A} i c_i$.  Either of these two upper bounds may be smaller, depending
on the situation given by $A, b$, and the $c_i$.

For a more complicated conditioning in which the $Z_i^*$ are not mutually
independent, consider random permutations on $n$ objects conditional on
having at least one cycle of length two or three. Here, $Z_2^*$ and $Z_3^*$
are dependent, although the {\it pair} $(Z_2^*,Z_3^*)$ and the variables
$Z_1^*,Z_4^*,Z_5^*,\ldots$ are mutually independent.  With $A=\{2,3\} 
\subseteq B
=\{1,2,\ldots,b\}$, we have $p=\B P(Z_2+Z_3 >0) = 1 - e^{-5/6}$ and
$d_b^*(n) \leq 3/(2p) F(n/b)$.  Thus, for example with $b=3$, 
the probability that a random permutation of $n$ objects 
is a derangement, given that
$C_2(n)+C_3(n)>0$, can be approximated by $\B P(Z_1^*=0)=1/e$, with
error at most $3/(2p) F(n/3)$.  Similarly, the probability that a 
random permutation of $n$ objects 
has a cycle of length 2, given that
$C_2(n)+C_3(n)>0$, can be approximated by $\B P(Z_2^*>0)=\B
P(Z_2>0|Z_2+Z_3>0) =(1-e^{-1/2})/(1-e^{-5/6})$, with
error again at most $3/(2p) F(n/3)$. 

The next example shows how to approximate easily the small component counts
for 
2--regular graphs by exploiting a decoupling result for the Ewens sampling
formula with parameter $\kappa = 1/2$.

\subsubsection{2-regular graphs}

The combinatorial structure known as `2--regular graphs' is the assembly
in which components are undirected cycles on three or more points, so
that
\eq\label{2reg-mi}
m_i = \frac{1}{2} (i-1)!\;\bone\{i \ge 3\}.
\en
Let $C_i^*(n)$ be the number of components of size $i$ in a random 
2--regular graph on $n$ points. A process that corresponds to this, with
the condition $\bone\{i \ge 3\}$ removed, is the Ewens sampling formula with
parameter $\kappa = 1/2$ described in Section
\ref{esfsect}. Observe that 
\[
\b C^*(n) \indist (\b C(n) | C_1(n) = C_2(n) = 0).
\]

The bound 
\[
d_{TV}( (C_1,\ldots,C_b), (Z_1,\ldots,Z_b)) \leq \frac {2 b}{n}
\]
is known from results of Arratia, Barbour and Tavar\'e (1992). We are
interested in how this translates into a bound on
\[
d_b^* \equiv d_{TV}( (C_3^*,\ldots,C_b^*), (Z_3,\ldots,Z_b)).
\]

With $A = \{1,2\}, B = \{1,2,\ldots,b\}$,
$d_A \leq 4/n, d_B \leq 2b/n, p = \B P(Z_1 = Z_2 = 0) = e^{-3/4}$,
the inequality in (\ref{tvbnd1}) guarantees that
\begin{eqnarray*}
d_b^* & \leq & \frac{1}{p} \left(\frac{d_A}{2}  + d_B\right) \\
&&\\
& \leq & e^{3/4} \left( \frac{2}{n} + \frac{2 b}{n} \right)\\
&&\\
& = & e^{3/4} \frac{2(b+1)}{n}. 
\end{eqnarray*}

For an example that shows that the conditioning event can
have probability tending to zero, consider 2--regular graphs conditioned
on having no cycles of size less than or equal to $t \equiv t(n) \geq 2$. The
previous example is the special case $t = 2$. For $b > t$, we have
\[
(C_{t+1}^*,\ldots,C_b^*) \indist (C_{t+1},\ldots,C_b | C_1 = \cdots =
C_t = 0).
\]
Now $d_A \leq 2t/n,\ d_B \leq 2b/n$, and
\[
p = \B P(Z_1= \cdots = Z_t=0) = \exp\left(-\frac{1}{2}(1+ \cdots+
1/t)\right)
\geq \frac{1}{\sqrt{e t}},
\]
so (\ref{tvbnd1}) establishes that 
\begin{eqnarray*}
d_b^* & \leq & \frac{1}{p}\left(\frac{d_A}{2} + d_B\right)\\
&&\\
& \leq & \sqrt{e t} \left(\frac{t}{n} + \frac{2b}{n}\right).
\end{eqnarray*}
This provides a useful bound provided that $\sqrt{t}b/n$ is small. Note
that both $t$ and $b$ may grow with $n$, as long as $t \leq b$.  
For example, conditional on no cycles of length less than or equal to $t
= \lfloor n^{2/3-\epsilon} \rfloor$ this approximation successfully
describes the distribution of the $k$ smallest cycles, for fixed $k$ as
$n \to \infty$, by using $b = n^{2/3}$. See Arratia and Tavar\'e (1992b,
Theorem 7) for related details.

 \section{Large deviation theory}\label{sect8}
\subsection{Biasing the combinatorial and independent processes}

A guiding principle of large deviation theory is that unlikely events of
the form $\{U \geq u \} $ or $\{ U \leq u \} $ or $ \{ U=u\}$, where the
target $u$ is far from $\B E U$, can be studied by changing the measure
$\B P$ to another measure $\B P_\theta$ defined by
\eq
\frac{d \B P_\theta}{d \B P} = \frac{\theta^U}{\B E \theta^U}.
\label{twist C}
\en
Observe that for $\theta=1$, the new measure $\B P_\theta$ coincides
with the original measure $\B P$, regardless of the choice of $U$. 
The parameter $\theta$ is chosen so that the average
value of $U$ under the new 
measure is $u$, i.e. $\B E_\theta U =u$. In  the literature on
large deviations and statistical mechanics (cf. Ellis, 1985), the
notation is usually 
$\theta \equiv e^\beta$, and our normalizing factor $\B E \theta^U$ is
expressed as the Laplace transform of the $\B P$-distribution of 
$U$, parameterized by $\beta$.

For the case of a combinatorial process $\b C(n) =(C_1(n),\ldots,C_n(n))$,
with  the total number of components 
$$
K \equiv K_n \equiv C_1(n) + \cdots + C_n(n)
$$ 
in the role of $U$, this says to change from the measure $\B P$, which
makes all possible structures equally likely, to the measure $\B
P_\theta$, which selects a structure with bias proportional to
$\theta^{\# {\rm components }}$.  The Ewens sampling formula
discussed in Section \ref{esfsect} is exactly this in the case of random
permutations, with $\kappa$ playing the role of $\theta$. This may 
easily be verified by comparing (\ref{esfdef}) to Cauchy's formula, the
special case $\kappa = 1$ of (\ref{esfdef}), in which the equality of
normalizing constants, with $\BE \kappa^{K_n} = \kappa_{(n)}$,
expresses a well known identity for Stirling numbers of the first kind.

Theorem  \ref{genequaldist} showed that many a combinatorial process is
equal in distribution to a process of independent random variables,
conditioned on the value of a weighted sum.  The next theorem asserts
that  this form is preserved
by the change of measure from large deviation theory, provided that $U$ is
also a weighted sum. 

As in the discussion before Theorem \ref{genequaldist}, the weight
function $\bu$, just like the weight function $\b w$, can take values in
$\BR$ or $\BR^d$.  In case the weights $\bu$, and hence the random
variable $U$, takes values in $\BR^d$ with $d>1$, we take $\theta >0$ to
mean that $\theta = (\theta_1,\ldots,\theta_d) \in (0,\infty)^d$, and
with $U=(U_1,\ldots,U_d), \theta^U$ represents the product
$\theta_1^{U_1}\cdots \theta_d^{U_d}$.

\begin{theorem}\label{twist dist thm}
Let $I$ be a finite set, and for $\alpha \in I$, let $C_\alpha$ and
$Z_\alpha$ be $\B Z_+$-valued random variables.
Let $\b w = (w(\alpha))_{\alpha \in
I}$ and $\b u = (u(\alpha))_{\alpha \in
I}$ be deterministic weight functions on $I$,
 with real values for $\b u$, let $T = \b w \cdot \b Z_I \equiv \sum_{\alpha \in
I} w(\alpha) Z_\alpha$, and let $U= \b u \cdot \b C_I$. Let $\B P$ be a
probability measure and $t$ be a
 constant such that, under $\B P$ the $Z_\alpha$ are mutually
independent, 
$\ \B P(T=t) >0$, and $\b C_I \indist (\b Z_I | T=t)$. Let $\theta >0$
be any constant such that the random variable  $Y \equiv
 \theta^{\b u \cdot \b Z_I}$ has $\BE Y
 < \infty$. Let $\B P_\theta$,
restricted to the sigma--field generated by $\b C_I$, be given by
(\ref{twist C}). Let $\BP_{\theta}$, restricted to the sigma-field
generated by $\bZ_I$, be given by 
$$
\frac{d\BP_{\theta}}{d\BP} = \frac{Y}{\BE Y},
$$
so that the  $Z_\alpha$ are  mutually independent under $\BP_{\theta}$ with 
\eq
\B P_\theta(Z_\alpha = k) = \frac{\theta^{u(\alpha)k}}{\B
E\theta^{u(\alpha) Z_\alpha}}\,\BP(Z_{\alpha} = k), k \geq 0.
\label{twist Z marginal}
\en
Then under $\B P_\theta, \ \ \b C_I \indist (\b Z_I |T=t)$, that is
\eq
\B P_\theta(\b C_I = \b a) = \B P_\theta(\b Z_I = \b a |T=t),
\en
for $\b a \in \B Z_+^I$.
\end{theorem}
\proof
For $\b a \in \B Z^I_+$,
\begin{eqnarray}
\B P_\theta (\b C_I = \b a) &=& (\B E \theta^U)^{-1} \ \ \theta^{\b u \cdot \b a}
 \ \ \B P(\b C_I = \b a)  \nonumber \\ && \nonumber \\
&=& (\B E \theta^U)^{-1}\ \  \theta^{\b u \cdot \b a}\ \
 \B P(\b Z_I = \b a| T=t)  \nonumber \\ && \nonumber \\
&=& (\B E \theta^U)^{-1} \ \B P(T=t)^{-1} \ \theta^{\b u \cdot \b a}
 \ \b 1(\b w \cdot \b a = t) \ 
\B P(\b Z_I = \b a)
\label{stept1}
\end{eqnarray}
Now
\[
\B P_\theta(\b Z_I=\b a) = \left(\B E \theta^{\b u \cdot \b
Z_I}\right)^{-1} \ 
\theta^{\b u \cdot \b a} \ \B P(\b Z_I = \b a)
\]
so that 
\eq
\B P_\theta(\b Z_I = \b a |T=t) = 
 \left(\B E \theta^{\b u \cdot \b Z_I}\right)^{-1} \ 
  \ \B P_\theta(T=t)^{-1} \ \theta^{\b u \cdot \b a}
\ \b 1(\b w \cdot \b a = t)  \ \B P(\b Z_I = \b a)
\label{stept2}.
\en
Comparing (\ref{stept1}) and (\ref{stept2}), we see both expressions are
probability densities on $\B Z_+^I$ which are 
proportional to the same function of $\b a$, and hence they are equal.
From this it also follows that the normalizing constants are equal,
which is written below with the combinatorial generating function on the
left, and the three factors determined by independent random variables
on the right:
\eq  
\B E \theta^U  =
 \B E \theta^{\b u \cdot \b Z_I}\ \frac{\B P_\theta(T=t)}{\B P(T=t)}.
\label{free and sexy}
\en
\hfill \qed

  For
the case $U = K_n$, the total number of components,
 the $\B P_\theta$ measure corresponds to the following generalization
of (\ref{Zassembly}) through (\ref{Zselection}).  For assemblies,
multisets, or selections, chosen with probability proportional to 
$\theta^{\# {\rm components} }, \b C(n) \indist ((Z_1,...,Z_n)|Z_1+2 Z_2 +\cdots
+ n Z_n = n)$ where the $Z_i$ are mutually independent. With $\theta, x
> 0$, for assemblies we have
\eq\label{twist Z dist}
 Z_i \mbox{ is }  \mbox{Poisson }
(\frac{m_i \ \theta \  x^i}{i!}),
\en
whereas for multisets we require $x \leq 1, \theta x < 1$ and then
$$
Z_i \mbox{ is negative binomial }(m_i,\theta \ x^i).
$$
Finally, for selections  
$$
Z_i \mbox{ is binomial }(m_i, \frac{\theta \  x^i}{1+\theta x^i}).
$$
In the general case, where $U = \bu \cdot \bC(n)$ is a weighted sum of
component counts, so that 
the selection bias is $\theta^{\b u \cdot \bC(n)}$, each factor 
$\theta$ in (\ref{twist Z dist}) above is replaced
by $\theta^{u(i)}$.  Furthermore, we observe that Theorems \ref{tvthm},
\ref{refinedthm}, and
 \ref{tvrefine} apply to $\B P_\theta$ in place of $\B P$.  For the
refinements in Section \ref{sect6}, for assemblies, multisets, and
selections respectively, the distribution of $Y_{ij}$ is Poisson
$(\theta^{u(i)} \ x^i/i!)$, Geometric ($\theta^{u(i)} \ x^i$),
 or Bernoulli ($\theta^{u(i)}
x^i/(1+\theta^{u(i)} x^i))$.

An example where such a bias is well known is the random graph model
${\mathcal G}_{n,p}$; Bollob\'as (1985). 
This corresponds to picking a labelled graph on $n$
vertices, where each of the potential edges is independently taken with
probability $p$; the unbiased case with all $2^{{n \choose 2}}$ graphs
equally likely is given by $p = 1/2$. We need something like the refined
setup of Section 6 to be able to keep track of components in terms of
the number of edges in addition to the number of vertices. Using the
full refinement of Section 6, $D_{ij}$ counts the number of components
on $i$ vertices having the $j$th possible structure, for $j = 1, \ldots,
m_i$, in some fixed enumeration of these. The weight function should be
$u(i,j)$ = \# edges in the $j$th possible structure on $i$ vertices. With
$\theta = p/(1-p)$, the $\BP_\theta$ law of ${\bf D}(n)$ is a
description of ${\mathcal G}_{n,p}$. A more natural refinement for this
example, intermediate between ${\bf C}$ and ${\bf D}$, would be the
process ${\bf A}$ with $A_{ik} = \sum_{j:u(i,j) = k} D_{ij}$, the number
of components with $i$ vertices and $k$ edges, for $k = i-1, \ldots, {i
\choose 2}$. As in (96) and (97), the total variation distances are
insensitive to the amount of refining. Presumably there are interesting
results about random graphs that could easily be deduced from estimates
of the total variation distance in Theorem 5.

One form of the general large deviation heuristic is that 
for a process $\b C$, conditioned on the event $\{ U \geq u \} $ 
where  $U$ is a functional of the process and $u > \B EU$, the $\BP-$law
of the conditioned process
is nicely approximated by the $\BP_\theta-$ law of $\b C$, where
$\theta$ is chosen so that $\B E_\theta U=u$.  We are interested in 
 the special case where the
functional $U$ is a weighted sum, and the distribution of $\b C$ 
under $\B P$
is that of an independent process $\b Z$ conditioned on the value of
another weighted sum $T$.  In this case,
 Theorem \ref{tvthm} yields a direct quantitative
handle on the quality of approximation by the $\B P_\theta$-distribution
of the independent process, provided we condition on 
the event $\{ U= u \} $ instead of the event $\{ U \geq u \} $.

\begin{theorem}\label{twist tvthm} 
Assume the hypotheses and notation of Theorems \ref{tvthm} and
 \ref{twist dist thm} combined.
For $B \subset I$ write $U_B \equiv \sum_{\alpha \in B} u(\alpha)
Z_\alpha$, so that $U_I \equiv \b u \cdot \b Z_I$.  
Write ${\mathcal L}_\theta$ for distributions governed by $\B
P_\theta$, so that the conclusion of Theorem \ref{twist dist thm} may be
written
\[
{\mathcal L}_\theta (\b C_I) = {\mathcal L}_\theta (\b Z_I | T=t),
\]
and Theorem \ref{tvthm} states that for $B \subset I$
\eq
d_{TV} ({\mathcal L}_{\theta}(\b C_B),{\mathcal L}_\theta (\b Z_B) \ ) =
d_{TV}({\mathcal L}_\theta (R_B|T=t), {\mathcal L}_\theta (R_B) \ ).
\en
 Assume that $u$ is such that $\B P(U=u)>0$. Then under the further
conditioning on $U=u$, 
$$
d_{TV} ({\mathcal L}_1(\b C_B|U=u),{\mathcal L}_\theta (\b Z_B) \ ) = \ \ \ \ \
\ \ \ 
$$
\eq
\ \ \ \ \ \
d_{TV}({\mathcal L}_\theta ((U_B,R_B)|U_I=u,T=t), {\mathcal L}_\theta ((U_B,R_B))).
\label{twisty}
\en
\end{theorem}
\proof
Observe first that 
\eq
\label{theta id}
{\mathcal L}_1(\b C_I|U=u) = {\mathcal L}_\theta(\b C_I|U=u),
\en
so that it suffices to prove (\ref{twisty}) with the subscript
$\theta$ appearing on all four distributions, i.e.
$$
d_{TV} ({\mathcal L}_\theta(\b C_B|U=u),{\mathcal L}_\theta (\b Z_B) \ ) =\ \ \
\ \ \ 
$$
\eq\label{twisty2}
\ \ \ \ \ \
d_{TV}({\mathcal L}_\theta((U_B,R_B)|U_I=u,T=t), {\mathcal L}_\theta ((U_B,R_B))).
\en
Observe next that this is a special case of Theorem \ref{tvthm}, but
with two--component weights $w^*(\alpha) \equiv (u(\alpha),w(\alpha))$
in the role of $w(\alpha)$. 
 For example, in the usual combinatorial
case, with $I =[n]$ and $w(i)=i$, and further specialized to $U=K_n=$
the total number of components, so that $u(i)=1$, we have that $\b w^*$
takes values in $\B R^2$, with $w^*(i)=(1,i)$. \hfill \qed

\medskip
\no{\bf Discussion.}  The proof of the previous theorem helps make it clear
that the free parameter $x$, such that ${\mathcal L}((Z_1,\ldots,Z_n)|T_n=n)$
does not vary with $x$, is analogous to the parameter $\theta$, such
that relation (\ref{theta id}) holds.  With this perspective, the
discussion of an appropriate choice of $x$ in Section \ref{sect4} and
Section \ref{sect5.2} is
simply giving details in some special cases of the general large deviation
heuristic.  Note that $T_n$ is a sufficient statistic for $x$, while $U$
is a sufficient statistic for $\theta$.

There are three distributions involved in the discussion above: the first
is ${\mathcal L}(\b C_I|U=u)$, corresponding to a combinatorial distribution
conditioned on the value of 
a weighted sum $U$, the second is ${\mathcal L}_\theta(\b C_I)$,
which is a biased version of the combinatorial distribution, and the 
third is ${\mathcal
L}_\theta(\b Z_I)$, which governs an independent process.  Theorem \ref{tvthm}, used
with Theorem \ref{twist dist thm}, compares the second and third of
these; Theorem \ref{twist tvthm}
above compares the first and third of these; and the following theorem
completes the triangle, by comparing the first and second distributions.

\begin{theorem}\label{leg3}
In the setup of Theorem \ref{twist tvthm},  for $B \subset I$,
\eq
d_{TV} ({\mathcal L}_1(\b C_B|U=u),{\mathcal L}_\theta (\b C_B) \ ) =
\label{leg3eq}
\en
\[ d_{TV}({\mathcal L}_\theta ((U_B,R_B)|U_I=u,T=t), {\mathcal L}_\theta
((U_B,R_B)|T=t \ )).
\]
\end{theorem}
\proof
By Theorem \ref{twist dist thm}, together with (\ref{theta id}),
 the left side of
(\ref{leg3eq}) is equal to $ d_{TV} ({\mathcal L}_\theta(\b Z_B|U_I=u,T=t),
{\mathcal L}_\theta (\b Z_B|T=t) \ ) $. We modify the second proof of Theorem
\ref{tvthm} as follows:  replace $\B P$  by $\B P_\theta$, use two--component 
weights, replace the original
conditioning $T=t$ by $U_I=u$, and then  further condition
 on $\{T=t \}$. Explicitly, the functional $h$ on $ \B Z_+^B$
defined by $h(\b a)= \sum_{\alpha \in B}
a(\alpha)(u(\alpha),w(\alpha))$ is a sufficient statistic, and the sign
of 
$\B P_\theta(\b Z_B= \b a|U_I=u,T=t) - \B P_\theta(\b Z_B=\b a | T = t) $ 
is equal to the sign of  $\B P_\theta((U_B,R_B)=h(\b a)|U_I=u,T=t)- 
\B P_\theta((U_B,R_B)=h(\b a)|T=t)$, i.e. the sign depends on $\b a$
only through the value of $h (\b a)$.
 \hfill \qed

Observe that Theorem \ref{twist tvthm} 
contains Theorem \ref{tvthm} as a special case, by
taking weights $u(\alpha) \equiv 0$ and target $u=0$, so that $\B
P_\theta = \B P$ and the extra conditioning event $\{ U=u \}$ has
probability one.  

\subsection{Heuristics for good approximation of conditioned
combinatorial structures}

The following applies to weighted sums $U$ in general, but to be concrete
we present the special case $U=K_n$.  
Let $K \equiv K_n$ be the total number of components of some assembly, multiset,
or selection of total weight $n$, and let some deterministic target $k \equiv
k(n)$ be given.  The goal is to describe an independent process to 
approximate $\b C(n)$, conditioned on
the event $\{K \geq k \}$, in case $k$ is large compared to $\B EK$; or
conditioned on 
the event $\{K \leq k \}$, in the opposite case; or more simply,
conditioned on the event $\{ K=k \}$. We accomplish this by picking the
free parameters $\theta$ and $x$ in (\ref{twist Z dist}) so that
simultaneously $\B E(Z_1+\cdots +Z_n)$ is close to $k$ and $\B ET_n$ is
close to $n$.  

For example, to study random permutations on $n$ objects, conditional on
having at least $5 \log n$ cycles, or conditional on having exactly
$\lfloor 5 \log n \rfloor$ cycles,  or conditional on having at most $0.3
\log n $ cycles, we propose using $x=1$, and $\theta = 5$ or $0.3$.  The
independent process with this choice of parameter should be a good
approximation for both the conditioned random permutations and for the
Ewens sampling formula.  As a corollary, the Ewens sampling formula
should be a good approximation for the conditioned permutations; see
Arratia, Barbour and Tavar\'e (1994).

For assemblies, multisets and selections in the logarithmic class
discussed in Section \ref{logsect}, in which $\BE Z_i \sim
\kappa/i$, biasing by $\theta^K$ yields $\BE_{\theta} Z_i \sim
\kappa \theta/i$, so that the Ewens sampling formula with
parameter $\kappa \theta$ is a useful approximation for the biased
measures. In particular, the heuristics (\ref{approx3b}) and
(\ref{approx4b}) should apply in the following form: for fixed 
$B \subseteq [n]$

 {\it In the case $\kappa \theta \neq 1$}
 \eq\label{approx3new}
 d_{TV}({\mathcal L}_{\theta}(\bC_B),{\mathcal L}_{\theta}(\bZ_B)) \sim 
\frac{1}{2} |\kappa \theta - 1|
 \frac{\BE_{\theta}|R_B - \BE_\theta R_B|}{n},
 \en
 
 {\it In the case $\kappa \theta = 1$}
 \eq\label{approx4new}
 d_{TV}({\mathcal L}_{\theta}(\bC_B),{\mathcal L}_{\theta}(\bZ_B)) = 
o\left(\frac{1}{n}\right).
 \en

For random permutations, for which $\kappa = 1$, with $B = \{1,2,\ldots,b\}$
the bound
$$
d_{TV}({\mathcal L}_{\theta}(\bC_B),{\mathcal L}_{\theta}(\bZ_B)) \leq c(\theta)
\frac{b}{n}
$$
was established via a particular coupling in Arratia, Barbour and
Tavar\'e (1992), and the asymptotic relation (\ref{approx3new}) has
been established by Arratia, Stark and Tavar\'e (1994).

To show how the parameters $x$ and $\theta$ may interact, we consider
random permutations with $k(n)$ further away from $\log n$. Assume that
$k(n)$ is given such that as $n \rightarrow \infty$,
\[
  k/ \log n \rightarrow \infty, \ \ \ \ k/n \rightarrow 0.
\]
Then we would take
\eq
\theta \equiv \theta(n) = \frac{k}{\log(n/k)}, \ \ \ \ \ x \equiv x(n) =
e^{-\theta /n}.
\label{ewen pair}
\en
Observe that $\theta/n \rightarrow 0$, so that $x \rightarrow 1$ and $1-
x \sim \theta/n$, and
$\theta \rightarrow \infty$, so that $x^n = \exp(-\theta) \rightarrow 0$.
Hence
\[ 
\B E T_n = \theta \sum_1^n x^i \sim \theta \sum_0^\infty x^i =
\theta \frac{1}{1-x} \sim n
\]
and 
\[
\B E K_n = \theta \sum_1^n \frac{x^i}{i} \sim -\theta \log(1-x) \sim
\theta \log (\frac{n}{\theta}) \sim k.
\]
With this choice of parameters $\theta$ and $x$ the independent Poisson
process $(Z_1,Z_2,\ldots)$ should be a good approximation for random
permutations, conditioned either on having exactly $k$ cycles, or on
having at least $k$ cycles.

 \section{The generating function connection and moments}\label{sect9}

In this section, we relate the probabilistic technique to the more
conventional one based on generating functions; Wilf (1990). 
One reason for this is to 
provide a simple method, based on an idea of Shepp and Lloyd (1966), for 
calculating moments of component counts for combinatorial structures. A 
second reason is to provide a framework within which detailed estimates and 
bounds for total variation distances can be obtained  by using  the results 
of Theorems \ref{tvthm} and \ref{twist tvthm}, together with analytic 
techniques such as Darboux's method or the transfer methods of
Flajolet and Odlyzko (1990).

  Throughout, we let 
$p(n,k)$ be the number of objects of weight $n$ having $k$ components,
so that
$p(n) = \sum_{k=1}^n p(n,k)$ is the number of objects of weight 
$n$. Finally, recall that $m_i$ is the number of available structures
for a component of size $i$.

\subsection{Assemblies}

We form the exponential generating functions
\eq\label{egf3}
\hat P(s,\theta) \equiv 1 + \sum_{n=1}^{\infty} \left(\sum_{k=1}^n p(n,k) \theta^k
 \right) \frac{s^n}{n!},
\en
\eq\label{egf1}
\hat P(s) \equiv 1 + \sum_{n=1}^{\infty} p(n) \frac{s^n}{n!} = \hat
P(s,1),
\en
and
\eq\label{egf2}
\hat M(s) \equiv  \sum_{n=1}^{\infty} m_n \frac{s^n}{n!}.
\en
For assemblies, (\ref{Nassembly}) gives
$$
p(n,k) = \sum_{\ba} N(n,\ba) = \sum_{\b a} n! \prod_{j=1}^n
\left(\frac{m_j}{j!}\right)^{a_j} \frac{1}{a_j!},
$$
where $\sum_{\b a}$ is over $\{\b a \in \B Z_+^n: \sum i a_i = n, \sum
a_i = k\}$. It follows that
\begin{eqnarray}
\hat P(s,\theta) & = & 1 + \sum_{n=1}^\infty \sum_{k=1}^n \sum_{\b a}
\prod_{j=1}^n \left( \frac{\theta m_j s^j}{j!}\right)^{a_j}
\frac{1}{a_j!} \nonumber \\
&&\nonumber \\
& = & \prod_{j=1}^{\infty} \exp\left(\frac{\theta m_j s^j}{j!}\right)
\nonumber \\
&& \nonumber \\
& = & \exp \left( \theta \hat M(s) \right) \label{egf4}.
\end{eqnarray}
Equation (\ref{egf4}) is the well-known exponential generating function
relation for assemblies (cf. Foata, 1974), which has as a special case the 
relationship
\eq\label{egf5}
\hat P(s) = \exp\left( \hat M(s) \right).
\en

Recall from Section 8 that in studying  large deviations  of $K_n$,
the number of components in the structure of total weight $n$, we were led 
to the measure $\BP_\theta$
corresponding to sampling with probability proportional to $\theta^{K_n}$. 
It follows from (\ref{Nassembly}) that there is a normalizing constant 
$p_\theta(n)$ such that
\begin{eqnarray*}
p_\theta(n)\,\B P_{\theta}(\b C(n) = \b a) & = &
 \theta^{a_1+ \cdots+a_n} N(n,\ba) \\
& = & n! x^{-n} \prod_{j=1}^n 
\left(\frac{\theta m_j x^j}{j!}\right)^{a_j} \frac{1}{a_j!}\ \bone\left(
\sum_{l=1}^n l a_l = n \right)
\end{eqnarray*}
for any $x > 0$. Clearly, 
\begin{eqnarray}
p_{\theta}(n) & = &  \sum_{k=1}^n p(n,k) \theta^k \nonumber \\
 & = & n!  [s^n] \hat P(s, \theta) \label{normconst1}\\
 & = & p(n) \B E(\theta^{K_n}) \label{normconst2},
\end{eqnarray}
where $\B E \equiv \BE_1$ denotes expectation with respect to the uniform measure $\BP \equiv
\BP_1$, corresponding to $\theta = 1$. 

Next we explore the connection with the
probability generating function (pgf) of the random variable $T_n \equiv
\sum_{j=1}^n j Z_j$, where the $Z_j$ are independent Poisson distributed
random variables with mean
$$
\BE_{\theta} Z_j \equiv \theta \lambda_j = \theta \frac{m_j x^j}{j!}.
$$
Recall  that the pgf of a Poisson-distributed random variable $Z$ with
mean $\lambda$ is
$$
\BE_{\theta} s^Z \equiv \sum_{j=0}^\infty \BP_{\theta}(Z=j) s^j \, = 
\exp(-\lambda (1-s)),
$$
so using the independence of the $Z_j$,
\begin{eqnarray*}
\BE_{\theta} s^{T_n} & = & \BE_{\theta} s^{\sum_{j=1}^n j Z_j} \\
&&\\
& = & \prod_{j=1}^n \BE_{\theta} \left(s^j\right)^{Z_j}\\
&&\\
& = & \exp\left( - \theta \sum_{j=1}^n \lambda_j (1 - s^j) \right)
\end{eqnarray*}
Thus
\begin{eqnarray*}
\BP_{\theta}(T_n = n) & = & [s^n]\, \BE_{\theta} s^{T_n} \\
&&\\
& = & \exp\left( -\theta \sum_{j=1}^n \lambda_j \right)\,[s^n] \exp
\left(\theta \sum_{j=1}^n \lambda_j s^j\right)\\
&&\\
& = & \exp\left( -\theta \sum_{j=1}^n \lambda_j \right)\,[s^n] \exp
\left(\theta \sum_{j=1}^\infty \lambda_j s^j\right)\\
&&\\
& = & \exp\left( -\theta \sum_{j=1}^n \lambda_j \right)\,[s^n] \exp\left( \theta \hat M(sx)\right)\\
&&\\
& = & \exp\left( -\theta \sum_{j=1}^n \lambda_j \right)\,[s^n] \hat P(sx, \theta),
\end{eqnarray*}
using (\ref{egf4}) at the last step. Thus, via  (\ref{normconst1}),
\eq\label{egf6}
\BP_{\theta}(T_n = n) = \exp\left(- \theta \sum_{j=1}^n \lambda_j\right) 
\frac{x^n p_{\theta}(n)}{n!},
\en
as can also be calculated from (\ref{assprobT=t}) and (\ref{free and sexy}) 
for the special case $U = K_n$.

The next result gives a simple expression for the joint moments 
of the component counts. We use the notation $y_{[n]}$ to denote the falling 
factorial $y(y-1)\cdots(y-n+1)$.
\begin{lemma}\label{moments ass}
For $(r_1,\ldots, r_b) \in \BZ_+^b$ with $m = r_1 + 2 r_2 +
 \cdots + b r_b$, we have
\eq\label{meanass}
\BE_{\theta} \prod_{j=1}^b (C_j(n))_{[r_j]} = \bone(m \leq n)\,x^{-m} 
\frac{n!}{p_{\theta}(n)} \frac{p_{\theta}(n-m)}{(n-m)!} \prod_{j=1}^b 
\left( \frac{\theta m_j x^j}{j!}\right)^{r_j}.
\en
\end{lemma}

\proof The key step is the substitution of $a_1,\ldots,a_b$ for $a_1-r_1,\ldots,a_b - r_b$ 
in the third equality below. For $m \leq n$, we have
\begin{eqnarray*}
\BE_{\theta} \prod_{j=1}^b (C_j(n))_{[r_j]} & = & \sum_{a_j \geq r_j,
j=1,\ldots,b}
\sum_{a_{b+1},\ldots,a_n: \sum j a_j = n} (a_1)_{[r_1]}\cdots (a_b)_{[r_b]}
\frac{n!}{x^n p_{\theta}(n)}\\ 
& & \ \ \times \prod_{j=1}^n \left( \frac{\theta m_j
x^j}{j!}\right)^{a_j} \frac{1}{a_j!} \\
&&\\
& = & \frac{n!}{x^n p_{\theta}(n)} \prod_{j=1}^b \left( \frac{\theta m_j
x^j}{j!}\right)^{r_j}\,\sum \sum \prod_{j=1}^b \left( \frac{\theta m_j
x^j}{j!}\right)^{a_j - r_j}\\ 
& & \ \ \times \frac{1}{(a_j - r_j)!} \prod_{j=b+1}^n 
\left( \frac{\theta m_j x^j}{j!}\right)^{a_j} \frac{1}{a_j!} \\
&&\\
& = & \frac{n!}{x^n p_{\theta}(n)} \prod_{j=1}^b \left( \frac{\theta m_j
x^j}{j!}\right)^{r_j}\,\sum_{a_1,\ldots,a_{n}: \sum j a_j = n-m} \\
&&\\
& & \ \ \times \prod_{j=1}^n 
\left( \frac{\theta m_j x^j}{j!}\right)^{a_j} \frac{1}{a_j!} \\
&&\\
& = & \frac{n!}{x^n p_{\theta}(n)} \prod_{j=1}^b \left( \frac{\theta m_j
x^j}{j!}\right)^{r_j}\, \frac{x^{n-m} p_{\theta}(n-m)}{(n-m)!}.
\end{eqnarray*}
\hfill \qed

\no{\bf Remark: } If $\{Z_i\}$ are  mutually independent Poisson random 
variables with $\BE_\theta Z_i = \theta m_i x^i / i!$, then the product term on the  right of equation (\ref{meanass}) is precisely 
$\BE_{\theta} \prod_{j=1}^b (Z_j)_{[r_j]}$.

\no{\bf Remark: } In the special case of permutations, in which $m_i =
(i-1)!$ and $p(n) = n!$, the normalizing constant $p_{\theta}(n)$  is
given by $p_{\theta}(n) = \theta(\theta+1)\cdots(\theta+n-1)$, and
equation (\ref{meanass}) reduces to 
$$
\BE_{\theta} \prod_{j=1}^b (C_j(n))_{[r_j]} = \bone(m \leq n)\, 
{{\theta+n-m-1} \choose {n-m}}{{\theta+n-1} \choose {n}}^{-1}
\prod_{j=1}^b \left( \frac{\theta}{j}\right)^{r_j},
$$
a result of Watterson (1974).

\subsection{Multisets}

For multisets, the (ordinary) generating functions are
\eq\label{ogf3}
P(s,\theta) \equiv 1 + \sum_{n=1}^{\infty} \left(\sum_{k=1}^n p(n,k) \theta^k
 \right) s^n,
\en
\eq\label{ogf1}
P(s) \equiv 1 + \sum_{n=1}^{\infty} p(n) s^n = P(s,1),
\en
and
\eq\label{ogf2}
M(s) \equiv  \sum_{n=1}^{\infty} m_n s^n.
\en
In this case, using (\ref{Nmultiset}) gives
$$
p(n,k) = \sum_{\ba} N(n,\ba) = \sum_{\b a} \prod_{j=1}^n
{{m_j+a_j-1} \choose {a_j}},
$$
the sum $\sum_{\b a}$ being over $\{\b a \in \B Z_+^n: \sum i a_i = n, \sum
a_i = k\}$. It follows that
\begin{eqnarray}
P(s,\theta) & = & 1 + \sum_{n=1}^\infty \sum_{k=1}^n \sum_{\b a}
\prod_{i=1}^n {{m_i+a_i-1} \choose {a_i}} (\theta s^i)^{a_i} \nonumber \\
&&\nonumber \\
& = & \prod_{i=1}^{\infty} (1- \theta s^i)^{-m_i}\label{ogf7}\\
&& \nonumber\\
& = & \exp\left(- \sum_{i=1}^{\infty} m_i \log(1 - \theta s^i)\right)
\nonumber\\
&&\nonumber \\
& = & \exp\left(
\sum_{i=1}^{\infty} m_i \sum_{j=1}^{\infty} \frac{ (\theta
s^i)^j}{j}\right)\nonumber \\
&&\nonumber \\
& = & \exp\left(
\sum_{j=1}^{\infty}\frac{\theta^j}{j} \sum_{i=1}^{\infty} m_i  s^{ij}\right)
\nonumber \\
&& \nonumber \\
& = & \exp\left(
\sum_{j=1}^{\infty}\frac{\theta^j}{j} M(s^j)\right). \label{ogf9}
\end{eqnarray}
See Foata (1974) and Flajolet and Soria (1990) for example.

Under the measure $\B P_{\theta}$, there is a normalizing constant $p_{\theta}(n)$ such that 
\begin{eqnarray*}
p_{\theta}(n)\,\B P_{\theta}(\b C(n) = \b a) & = & 
\prod_{i=1}^n {{m_i +
a_i -1} \choose {a_i}} \theta^{a_i}\,\bone\left(\sum_{l=1}^n l a_l =
n\right) \\
& = & x^{-n} \prod_{l=1}^n (1 - \theta x^l)^{-m_l}  
\prod_{i=1}^n {{m_i +
a_i -1} \choose {a_i}} (1-\theta x^i)^{m_i} (\theta x^i)^{a_i}\\
& & \ \ \times 
\bone\left(\sum_{l=1}^n l a_l = n\right),
\end{eqnarray*}
for any $0 < x < 1$. Indeed,
\eq\label{ogf11}
p_{\theta}(n) = p(n) \B E_1(\theta^{K_n}) = [s^n] P(s,\theta),
\en
where $p_{\theta}(0) \equiv 1.$

 In this case, the relevant $Z_j$ are independent negative binomial random variables with parameters $m_i$ and $\theta x^i$ and pgf
$$
\BE_{\theta} s^{Z_i} = \left(\frac{1-\theta x^i}{1-\theta x^is}\right)^{m_i}.
$$
Using the independence of the $Z_j$ once more, the pgf of 
$T_n$ may be found as
\begin{eqnarray}
\BE_{\theta} s^{T_n} & = & \prod_{i=1}^n \BE_{\theta} \left(s^i\right)^{Z_i}
\nonumber \\
&& \nonumber \\
& = & \prod_{i=1}^n \left(\frac{1-\theta x^i}{1-\theta
(xs)^i}\right)^{m_i}\nonumber \\
&& \nonumber \\
& = & \left(\prod_{i=1}^n (1-\theta x^i)^{m_i}\right) 
\prod_{i=1}^n (1-\theta (xs)^i)^{-m_i} \label{ogf13}.
\end{eqnarray}
Using (\ref{ogf7}), we see that
\begin{eqnarray*}
\BP_{\theta}(T_n = n) & = & [s^n]\, \BE_{\theta} s^{T_n} \\
&&\\
& = & \left(\prod_{i=1}^n (1-\theta x^i)^{m_i}\right) [s^n] \exp\left(-
\sum_{i=1}^n m_i \log(1 - \theta (xs)^i)\right)\\
&&\\
& = & \left(\prod_{i=1}^n (1-\theta x^i)^{m_i}\right) [s^n] \exp\left(-
\sum_{i=1}^{\infty} m_i \log(1 - \theta (xs)^i)\right)\\
&&\\
& = & \left(\prod_{i=1}^n (1-\theta x^i)^{m_i}\right) [s^n]
P(sx,\theta),
\end{eqnarray*}
so that
\eq\label{ogf12}
\BP_{\theta}(T_n = n) = \left(\prod_{i=1}^n (1-\theta x^i)^{m_i}\right) 
x^n p_{\theta}(n).
\en
Equation (\ref{ogf12}) can also be calculated from (\ref{assprobT=t}) 
and (\ref{free and sexy}) for the special case $U = K_n$.

In order to calculate moments of the component counts $\b C(n)$, it is
convenient to use a variant on a theme of Shepp and Lloyd (1966). We
assume that $M(s)$ has positive radius of convergence, $R$.
As above, let $Z_1,Z_2,\ldots$ be mutually independent negative binomial random
variables, $Z_i$ having parameters $m_i$ and $\theta x^i$, where $0 < x
< \min\{R,1,\theta^{-1}\}$. Let $T_\infty \equiv \sum_{i=1}^{\infty} i Z_i$. 
Note that $T_{\infty}$ is almost surely finite, because
$$
\B E_{\theta} T_{\infty} = \sum_{i=1}^{\infty} \frac{i m_i \theta x^i}{1
- \theta x^i} \leq \frac{\theta x}{1 - \theta x} M^{\prime}(x) < \infty.
$$
The distribution of $T_{\infty}$ follows from (\ref{ogf7}), (\ref{ogf13}) and 
monotone convergence since 
$$
\BE_{\theta} s^{T_{\infty}} = \frac{P(sx,\theta)}{P(x,\theta)}.
$$
Hence
\eq\label{ogf14}
\BP_{\theta}(T_{\infty} = n) = x^n p_{\theta}(n) / P(x,\theta),
n=0,1,\ldots.
\en
Further, for $\b a \in \BZ_+^n$ and  $\b Z(n) \equiv (Z_1,\ldots,Z_n)$
\eq\label{ogf15}
\BP_{\theta}(\b C(n) = \b a) = \BP_{\theta}(\b Z(n) = \b a | T_{\infty} =
n).
\en
This follows from the statement (\ref{twist Z dist}) that
$$
\BP_{\theta}(\bC(n) = \ba) = \BP_{\theta}(\bZ(n) = \ba | T_n = n),
$$
and the observation that
\begin{eqnarray*}
\BP_{\theta}(\bZ(n) = \ba | T_n = n) & = & \frac{\BP_{\theta}(\bZ(n) = 
\ba, T_n = n)}{\BP_{\theta}(T_n = n)} \\
& = & \frac{\BP_{\theta}(\bZ(n) = \ba, T_n = n) \BP_{\theta}(Z_{n+1} = 
Z_{n+2} = \cdots = 0)}{ \BP_{\theta}(T_n = n) \BP_{\theta}(Z_{n+1} = Z_{n+2} = \cdots = 0)} \\
& = & \frac{\BP_{\theta}(\bZ(n) = \ba, T_{\infty} = n)}{\BP_{\theta}
(T_{\infty} = n)} \\
& = & \BP_{\theta}(\bZ(n) = \ba | T_{\infty} = n).
\end{eqnarray*}

Now let $\Phi: \BZ_+^{\infty} \to \BR$, and set $\bC_n \equiv
(C_1(n),\ldots,C_n(n),0,0,\ldots)$ with $\bC_0 \equiv (0,0,\ldots)$. 
The aim is to find an easy way to
compute $\BE_{\theta}^n( \Phi) = \BE_{\theta}\Phi(\bC_n)$.
It is convenient to use the notation $\BE_{x,\theta}$ to denote
expectations computed under the independent negative binomial measure
with parameters $x$ and $\theta$. Shepp and Lloyd's method in the
present context is the observation, based on (\ref{ogf14}) and
(\ref{ogf15}), that $\BE_{x,\theta}(\Phi | T_{\infty} = n) = \BE_{\theta}^n(\Phi)$, so that
\begin{eqnarray}
\BE_{x,\theta} (\Phi) & = & \sum_{n=0}^\infty \BE_{x,\theta}( \Phi |
T_{\infty} = n) \BP_{\theta}(T_{\infty} = n) \nonumber \\
&& \nonumber \\
& = & \sum_{n=0}^\infty \BE^n_{\theta}( \Phi)
x^n p_{\theta}(n)/ P(x,\theta). \label{sandl1}
\end{eqnarray}
This leads to the result that
\eq\label{ogf16}
\BE^n_{\theta}(\Phi) = \frac{[x^n] \BE_{x,\theta}(\Phi)
P(x,\theta)}{p_{\theta}(n)}.
\en

For $r \geq 1, jr \leq n$, we use this method to calculate
the falling factorial moments 
$\BE_{\theta} (C_j(n))_{[r]}$.  This  determines all moments, 
since $C_j(n)_{[r]} \equiv 0$ if $jr > n$.
In this case $\Phi(x_1,x_2,\ldots) = (x_j)_{[r]}$, and
\begin{eqnarray*}
\BE_{x,\theta}(\Phi) & = & \BE_{x,\theta} (Z_j)_{[r]} \\
& = & \frac{\Gamma(m_j+r)}{\Gamma(m_j)} \left( \frac{\theta x^j}{1 -
\theta x^j}\right)^r.
\end{eqnarray*}
Hence we have
\begin{eqnarray}
\BE_{\theta}^n(\Phi) & = & \frac{\Gamma(m_j+r)}{p_{\theta}(n)
\Gamma(m_j)} [x^n] P(x,\theta) \left( \frac{\theta x^j}{1 -
\theta x^j}\right)^r \nonumber\\
&&\nonumber\\
& = & \frac{\theta^r \Gamma(m_j+r)}{p_{\theta}(n) \Gamma(m_j)} [x^{n-
rj}] P(x,\theta) \sum_{l=0}^{\infty} {{r+l-1} \choose {l}} \theta^l
x^{jl} \nonumber\\
&& \nonumber\\
& = & \frac{\theta^r \Gamma(m_j+r)}{p_{\theta}(n) \Gamma(m_j)}
\sum_{l=0}^{\lfloor n/j \rfloor- r} {{r+l-1} \choose {l}} \theta^l
p_{\theta}(n - jr - jl) \nonumber \\
&& \nonumber\\
& = & \frac{\Gamma(m_j+r)}{p_{\theta}(n) \Gamma(m_j)}
\sum_{m = r}^{\lfloor n/j \rfloor} {{m-1} \choose {r-1}} \theta^m
p_{\theta}(n - jm).
\end{eqnarray}

\no{\bf Remark: } See Hansen (1993) for related material. The Shepp and 
Lloyd method can also be used in
the context of assemblies, for which (\ref{ogf14}) holds with
\eq\label{ogf26}
\BP_{\theta}(T_{\infty} = n) = \frac{x^n}{n!} p_{\theta}(n) / \hat P(x,\theta),
n=0,1,\ldots\ .
\en
This provides another proof of Lemma \ref{moments ass}.
See Hansen (1989) for the case of random mappings, and Hansen (1990) for
the case of the Ewens sampling formula.

\subsection{Selections}

The details for the case of selections are similar to those for
multisets. Most follow by replacing $\theta$ and $m_i$ by $-\theta$ and
$-m_i$ respectively in the formulas for multisets. First, we have from (\ref{Nselection})
$$
p(n,k) = \sum_{\ba} N(n, \ba) = \sum_{\b a} \prod_{j=1}^n
{{m_j} \choose {a_j}},
$$
the sum $\sum_{\b a}$ being over $\{\b a \in \B Z_+^n: \sum i a_i = n, \sum
a_i = k\}$. Therefore 
\begin{eqnarray}
P(s,\theta) & = & 1 + \sum_{n=1}^\infty \sum_{k=1}^n \sum_{\b a}
\prod_{i=1}^n {{m_i} \choose {a_i}} (\theta s^i)^{a_i} \nonumber \\
&&\nonumber \\
& = & \prod_{i=1}^{\infty} (1 + \theta s^i)^{m_i}\label{sgf7}\\
&& \nonumber\\
& = & \exp\left(
\sum_{j=1}^{\infty}\frac{(-1)^{j-1} \theta^j}{j} M(s^j)\right), \label{sgf9}
\end{eqnarray}
the last following just as (\ref{ogf9}) followed from (\ref{ogf7}).
 See Foata (1974), Flajolet and 
Soria (1990) for example.

Under the measure $\B P_{\theta}$, there is a normalizing constant $p_\theta(n)$  such that 
\begin{eqnarray*}
p_{\theta}(n)\,\B P(\bC(n) = \b a) & = & \prod_{i=1}^n {{m_i}
\choose {a_i}} \theta^{a_i}\,\bone\left(\sum_{l=1}^n l a_l =
n\right) \\
& = & x^{-n} \prod_{l=1}^n (1 + \theta x^l)^{m_l}  
\prod_{i=1}^n {{m_i} \choose {a_i}} (1+ \theta x^i)^{-m_i} (\theta x^i)^{a_i}\\
& & \ \ \times 
\bone\left(\sum_{l=1}^n l a_l = n\right),
\end{eqnarray*}
for any $ x>0$;
$p_{\theta}(n)$ is given in (\ref{ogf11}) once more.
 In this case, the $Z_j$ are independent 
binomial random variables with pgf
\eq\label{sgf20}
\BE_{\theta} s^{Z_i} = \left(\frac{1+ \theta x^i s}{1+ \theta x^i}\right)^{m_i},
\en
and the pgf of $T_n$ is
\begin{eqnarray}
\BE_{\theta} s^{T_n} & = & 
\left(\prod_{i=1}^n (1+\theta x^i)^{-m_i}\right) 
\prod_{i=1}^n (1+\theta (xs)^i)^{m_i} \label{sgf13}.
\end{eqnarray}
It follows from (\ref{ogf7}) that
\begin{eqnarray*}
\BP_{\theta}(T_n = n) & = & 
\left(\prod_{i=1}^n (1+ \theta x^i)^{- m_i}\right) [s^n]
P(sx,\theta),
\end{eqnarray*}
so that
\eq\label{sgf12}
\BP_{\theta}(T_n = n) = \left(\prod_{i=1}^n (1+\theta x^i)^{-m_i}\right) 
x^n p_{\theta}(n).
\en

The joint moments of the counts can be calculated using the Shepp and
Lloyd construction once more. In particular, equations (\ref{ogf14}) and
(\ref{ogf15}) hold, and we can apply (\ref{ogf16}) with
$\BE_{x,\theta}(\Phi)$ denoting expectation with respect to independent
binomial random variables $Z_1,Z_2,\ldots$ with distribution determined
by (\ref{sgf20}). 

As an example, we use this method to calculate
$\BE_{\theta} (C_j(n))_{[r]}$ for $r \geq 1, jr \leq n$. Since
\begin{eqnarray*}
\BE_{x,\theta}(\Phi) & = & \BE_{x,\theta} (Z_j)_{[r]} \\
& = & (m_j)_{[r]} \left( \frac{\theta x^j}{1 + \theta x^j}\right)^r,
\end{eqnarray*}
from (\ref{ogf16}) we have
\begin{eqnarray}
\BE_{\theta}^n(\Phi) & = & \frac{(m_j)_{[r]}}{p_{\theta}(n)}
[x^n] P(x,\theta) \left( \frac{\theta x^j}{1 +
\theta x^j}\right)^r \nonumber\\
&&\nonumber\\
& = & \frac{(m_j)_{[r]}}{p_{\theta}(n)}
\sum_{m = r}^{\lfloor n/j \rfloor} {{m-1} \choose {r-1}} (-1)^{m-r} \theta^m
p_{\theta}(n - jm).
\end{eqnarray}

\subsection{Recurrence relations and numerical methods}

We saw in Theorems \ref{tvthm} and \ref{twist tvthm} that for 
any $B \subseteq [n]$, the total
variation distance between $\bC_B$ and $\bZ_B$ can be expressed in terms
of the distributions of random variables $S_B$ and $R_B$ defined by
\eq\label{egf8}
S_B = \sum_{i \in [n] - B} i Z_i,
\en
and
\eq\label{egf7}
R_B = \sum_{i \in B} i Z_i \equiv S_{[n]-B}.
\en
Specifically,
\begin{eqnarray}
d_{TV}({\mathcal L}_{\theta}(\bC_B), {\mathcal L}_{\theta}(\b Z_B)) & = &
\frac{1}{2} \B P_{\theta}(R_B>n) \label{tvcalc2} \\
& &\ \  + \frac{1}{2} \sum_{r=0}^n \B P_{\theta}( R_B=r) 
 \left| \frac{\B P_{\theta}( S_B=n-r)}{\B P_{\theta}(T_n=n)} - 1 \right|. 
\nonumber
\end{eqnarray}

A direct attack on estimation of $d_{TV}({\mathcal L}_{\theta}(\bC_B), {\mathcal L}_{\theta}(\bZ_B))$  can be based
on a generating function approach to the asymptotics (for large $n$) of
the terms in (\ref{tvcalc2}). In the setting of assemblies, this uses
the result before (\ref{egf6})
for $\BP_{\theta}(T_n = n)$, and the fact that for $k \geq 0$
\begin{eqnarray}
\BP_{\theta}(S_B = n-k) & = & [s^{n-k}]\,\exp \left(- \theta \sum_{i \in [n]-B}
\lambda_i(1 - s^i)\right) \nonumber\\
&&\nonumber\\
& = & \exp\left(- \theta \sum_{i \in [n]-B} \lambda_i\right)\,[s^{n-k}]
\exp \left( \theta \sum_{i \in [n]-B} \lambda_i s^i + \theta 
\sum_{i > n} \lambda_i s^i\right) \nonumber \\
&&\nonumber\\
& = & \exp\left(- \theta \sum_{i \in [n]-B} \lambda_i\right)\,[s^{n-k}]
\hat P(sx,\theta)\, \exp \left( - \theta \sum_{i \in B} \lambda_i
s^i\right). \nonumber\\
&&\label{recur0}
\end{eqnarray}
For applications of this technique, see Arratia, Stark and Tavar\'e
(1994), and Stark (1994b).

It is also useful to have a recursive method for calculating the
distribution of $R_B$ for any $B \subseteq [n]$. For assemblies,
\eq\label{recur1}
\BE_{\theta} s^{R_B} = \exp\left(- \sum_{i \in B} \theta \lambda_i
\right)\, \exp
\left( \sum_{i \in B} \theta \lambda_i s^i\right).
\en
Write 
$$
G_B(s) = \sum_{i \in B} \theta \lambda_i,
$$
and
$$
F_B(s) = \exp G_B(s) \equiv \sum_{k=0}^\infty q_B(k) s^k,
$$
with $q_B(0) \equiv 1$. Differentiating with respect to $s$ shows that $s
F^\prime_B(s) = s G^\prime_B(s)\,F_B(s)$ (cf. Pourahmadi, 1984), and
equating coefficients of $s^k$ gives 
$$
k q_B(k) = \sum_{i=1}^k g_B(i) q_B(k-i),\ k=1,2,\ldots
$$
where
\eq\label{recur2}
g_B(i) = \theta i \lambda_i \,\bone(i \in B).
\en
Since $p_B(k) \equiv \BP_{\theta}(R_B = k) = \exp(- G_B(1)) q_B(k)$,
we find that 
\eq\label{recur3}
k p_B(k) = \sum_{i=1}^k g_B(i) p_B(k-i),\ k=1,2,\ldots
\en
with $p_B(0) = \exp(-G_B(1))$.
The relation (\ref{recur3}) has been exploited in the case of uniform
random permutations ($\theta = 1, \lambda_i = 1/i$) by Arratia and Tavar\'e
(1992a).

For multisets, the analog of (\ref{recur0}) is
\begin{eqnarray}
\BP_{\theta}(S_B = n-k) & = & \left(\prod_{i \in [n]-B} (1- \theta 
x^i)^{m_i}\right) 
[s^{n-k}] \prod_{i \in [n]-B} (1-\theta (xs)^i)^{-m_i}\nonumber\\
&&\nonumber\\
& = & \left(\prod_{i \in [n]-B} (1-\theta x^i)^{m_i}\right) [s^{n-k}]
\prod_{i \in [n]-B} (1-\theta (xs)^i)^{-m_i} \nonumber\\
&&\nonumber\\
& & \ \ \times  \prod_{i>n} (1-\theta (xs)^i)^{-m_i}\nonumber\\
&&\nonumber\\
& = & \left(\prod_{i \in [n]-B} (1- \theta x^i)^{m_i}\right) [s^{n-k}]
P(sx, \theta) \prod_{i \in B} (1-\theta (xs)^i)^{m_i}.\nonumber\\
&&\label{recur4}
\end{eqnarray}

To develop a recursion for $p_B(k) \equiv \BP_{\theta}(R_B = k)$, we can
use logarithmic differentiation; cf Apostol (1976), Theorem 14.8.
First, we have
\eq\label{recur5}
\BE s^{R_B} = \prod_{i \in B} (1 - \theta  x^i)^{m_i}\,
\prod_{i \in B} (1 -  \theta  x^i s^i)^{- m_i}.
\en
Define 
$$
G_B(s) = \sum_{i \in B} m_i s^i,
$$
and
$$
F_B(s) = \prod_{i \in B} (1 -  \theta  x^i s^i)^{- m_i} 
\equiv \sum_{k=0}^\infty q_B(k) s^k,
$$
with $q_B(0) = 1$. Then
$$
\log F_B(s) = \sum_{j=1}^{\infty} \frac{\theta^j}{j} G_B( (xs)^j).
$$
Differentiating with respect to $s$ and simplifying shows that 
$$
s F^\prime_B(s) = \left( \sum_{i \geq 1} g_B(i) s^i \right)\,F_B(s),
$$
where
\eq\label{recur6}
g_B(i) = x^i \sum_{k|i} k m_k \theta^{i/k}\,\bone(k \in B).
\en
Equating coefficients of $s^k$ gives 
$$
k q_B(k) = \sum_{i=1}^k g_B(i) q_B(k-i),\ k=1,2,\ldots.
$$
Since $p_B(k) \equiv \BP_{\theta}(R_B = k) = \prod_{i \in B} (1 - \theta
x^i)^{m_i}\, q_B(k)$,
it follows that 
\eq\label{recur7}
k p_B(k) = \sum_{i=1}^k g_B(i) p_B(k-i),\ k=1,2,\ldots
\en
with $p_B(0) = \prod_{i \in B} (1 - \theta x^i)^{m_i}$.

For selections, we have the following identity, valid for $k \geq 0$:
$$
\BP_{\theta}(S_B = n-k)  = 
\left(\prod_{i \in [n]-B} (1+ \theta x^i)^{-m_i}\right) [s^{n-k}]
P(sx, \theta) \prod_{i \in B} (1+\theta (xs)^i)^{-m_i}.
$$
If we define $p_B(k) \equiv \BP_\theta(R_B =k)$, then from equation
(\ref{recur7}) we obtain 
\eq\label{stnew}
k p_B(k) = \sum_{i=1}^k g_B(i) p_B(k-i),\ k=1,2,\ldots
\en
where 
$$
g_B(i) = - x^i \sum_{k|i} k m_k (- \theta)^{i/k} \bone(k \in B),
$$
and 
$$
p_B(0) = \prod_{i \in B} (1 + \theta x^i)^{- m_i}.
$$

 \section{Proofs by overpowering the conditioning}\label{sect10}

The basic strategy for making the relation $\b C_I \indist (\b Z_I|T=t)$
into a useful approximation is to pick the free parameter $x$ in the
distribution of $\bZ_I$ so that the conditioning is not severe, i.e. so
that $\B P(T=t)$ is not too small.  It is sometimes possible to get useful upper
bounds on events involving the combinatorial process $\b C_I$ by
combining upper bounds on the probability of the same event for the
independent process, together with lower bounds for $\B P(T=t)$.
The formal description of this strategy is given by the following lemma.

\begin{lemma}\label{dumbdumb}
Assume that  $\b C_I \indist (\b Z_I|T=t)$
and that $h$ is a nonnegative functional of these processes, i.e.
\[
h: \B Z_+^I \rightarrow \B R_+.
\]
Then
\eq
\label{dumb but useful}
\B E h(\b C_I) \leq \frac{\B E h(\b Z_I)}{\B P(T=t)}.
\en
\end{lemma}
\proof
\[
\B E h(\b C_I) =\frac{\B E(h(\b Z_I) \bone (T=t))}{\B P(T=t)} \leq
  \frac{\B E h(\b Z_I)}{\B P(T=t)}.
\]
\hfill \qed

\subsection{Example: partitions of a set}
Recall that partitions of a set is the assembly with $m_i = 1$ for all
$i$.  
Following the discussion in Section \ref{sect5.3}
we take $x \equiv x(n)
=\log n - \log \log n + \cdots$ to be the solution on $x
e^x=n$, so that for
$i=1,2,\ldots,n$, $Z_i$ is Poisson distributed, with mean and variance
$\lambda_i = x^i/i!.$  
With this choice of $x$, we have
\[
\B E T_n = \sum_1^n i \lambda_i \sim x e^x  = n
\]
and
\eq\label{Moser2}
\sigma_n^2 \equiv {\rm var}(T_n) = \sum_1^n i^2 \lambda_i \sim n \log n.
\en
By combining (\ref{assprobT=t}) with the asymptotics for the Bell numbers
given in Moser and Wyman (1955), and simplifying, we see that
\eq \label{Moser}
\B P(T_n = n) \sim \frac{1}{\sqrt{2 \pi n \log n}} \sim \frac{1}{\sqrt{2
\pi } \ \sigma_n},
\en
which is easy to remember, since it agrees with what one would guess
from the local central limit heuristic.

Write $U_n =Z_1+Z_2+\cdots+Z_n$, so that 
the total number of blocks $K_n$ satisfies $K_n \indist (U_n|T_n=n)$.
Harper (1967) proved that $K_n$ is
asymptotically normal with mean $n/x$ and variance $n/x^2$.  We observe
that this contrasts with the unconditional behavior: $U_n$ is
asymptotically normal with mean $n/x$, like $K_n$, and variance $n/x$,
unlike $K_n$.  Since $U_n$ is
Poisson, it has equal mean and variance. Harper's result says that
conditioning on $T_n=n$ reduces the variance of $U_n$ by a factor
asymptotic to $\log n$. 

Note that $Z_1$ is Poisson with parameter $x \sim \log n$, and hence
the distribution of $Z_1$ is asymptotically normal with mean and
variance $\log n$.  Note also that 
the Poisson parameters $\lambda_i = x^i/i!$ are themselves
proportional to $\B P(Z_1=i)$; in fact for $i \geq 1$
\[
\lambda_i = e^x \B P(Z_1 =i) = \frac{n}{x} \B P(Z_1=i).
\]
We can use the normal
approximation for $Z_1$ to see that, for fixed $a<b$, as $n \rightarrow
\infty$, 
\[
\sum_{a \sqrt{\log n} < i-\log n < b \sqrt{\log n}} \lambda_i \ \sim \  
  \frac{n}{\log n} \ \int_a^b \frac{1}{\sqrt{2 \pi}} e^{-u^2/2} du.
\]  
Informally, the relatively 
large values of $\lambda_i$ occur for $i$ within a few
$\sqrt{\log n}$ of $\log n$.

\subsubsection{The size of a randomly selected block}
A result similar to the following appears  as Corollary 3.3 in DeLaurentis and Pittel
(1983).
The size $D \equiv D_n$ of ``a randomly selected component'' of a
random assembly on $n$ elements is defined by a two step procedure:
first pick a random assembly, then pick one of its $K_n$ components,
each with probability $1/K_n$.  The same definition applies to the case
of random multisets or selections of weight $n$.  

Given $1 \leq b \leq n$, consider the functional $h:\B Z_+^n \rightarrow
[0,1]$ defined by
\[
h(\b a) = \left( \sum_{i \leq b} a_i \right)  \ / 
 \left( \sum_{i \leq n} a_i \right) ,
\]
with $h(0,0,\ldots,0)$ defined to be 1. 
The distribution of the size of a randomly selected component is
determined by
\[
\B P(D_n \leq b) \ = \ \B E h(\b C(n)).
\]
Define $U_b = Z_1 + \cdots + Z_b$, so that $h(Z_1,\ldots, Z_n) =
U_b/U_n$ and 
\[
\B P(D_n \leq b) \ = \ \B E h((Z_1,\ldots,Z_n)|T_n=n) = \B E \left(
\left. \frac{U_b}{U_n} \right| T_n=n \right) .
\]
Let $\epsilon >0$ and $\rho >1$ be given.  Let $1 \leq b \leq n$ such that
\eq \label{qrelation}
q \equiv \B P(Z_1 \leq b) \ \in [2 \epsilon,1-2 \epsilon].
\en
Now for all $n \geq n(\epsilon,\rho)$ we have $\B E U_b > \epsilon \ n
/\log n$ and $\B E U_b / \B E U_n \ \in [q/\rho,q \rho]$.
Large deviation theory says that for $\rho >1$ there is a constant
$c=c(\rho) >0$ such that if $Y$ is Poisson with parameter $\lambda$,
then $\B P(Y/\lambda \leq 1/\rho) \leq \exp(-\lambda c)$ and $\B
P(Y/\lambda \geq \rho) \leq \exp(-\lambda c)$. (In fact, the optimal $c$
is given by $c(\rho)= \min(1+\rho \log \rho - \rho, 1- \rho^{-1} \log
\rho - \rho^{-1})$, with the two terms in the minimum corresponding
respectively  to large
deviations above the mean and below the mean.) Putting
these together, using the large deviation bounds once with $U_b$ as $Y$
and a second time with $U_n$ as $Y$, we have for $n \geq n(\epsilon,\rho)$
\[
\B P(\frac{U_b}{U_n} \notin [q/\rho^3,q \rho^3] \ ) \leq
 2 \exp( -c(\rho) \epsilon \ n / \log n).
\]
Since the functional $h$ takes values in $[0,1]$, this proves, for $n
\geq n(\epsilon,\rho)$,
\eq  \label{with a rate}
\left| \B P(D_n \leq b) - q \right| \  \leq q(\rho^3 -1) +
 2 \exp( -c(\rho) \epsilon \ n / \log n)/\B P(T_n=n).
\en 
In terms of Lemma \ref{dumbdumb}, the above argument involves the
functional $h^*$ defined by $h^*(\b a) = \bone (h(\b a) \notin [q/\rho^3, q
\rho^3])$. The inequality (\ref{with a rate}) not only proves that $D_n$
is asymptotically normal with mean and variance $\log n$, but also provides
an upper bound on the Prohorov distance between the distributions
of $D_n$ and  $Z_1$.

\subsubsection{The size of the block containing a given element}
In the case of
assemblies, it is possible that someone describing ``a randomly selected
component'' has in mind the component containing a randomly selected
element, where the element and the assembly are chosen independently. 
This includes, for example, the case where the element is
deterministically chosen, say it is always 1.  Let  $D^*_n$ be the size
of the component containing 1, in a random assembly on the set $[n]$. 

The two notions of ``a
randomly selected component'' can be very far apart.  For example, with
random permutations, $D^*_n$ is uniformly distributed over
$\{1,2,\ldots,n\}$, while the size $D_n$ of a randomly selected cycle is
such that $\log D_n / \log n$ is approximately uniformly distributed
over $[0,1]$.  For random partitions of a set, the argument below proves
that $D_n$ and $D_n^*$ are close in distribution, because both
distributions are close to Poisson with parameter $x$, where $xe^x=n$.

Given $1 \leq b \leq n$, consider the functional $g:\B Z_+^n \rightarrow
[0,1]$ defined by
\[
g(\b a) = \frac{1}{n} \sum_{i \leq b}i a_i  .
\]
The distribution of the size of the component containing a given element
is
determined by
\[
\B P(D^*_n \leq b) \ = \ \B E g(\b C(n)).
\]
Define $R_b = Z_1 + 2 Z_2 + \cdots +b Z_b$, so that $g(Z_1,\ldots, Z_n) =
R_b/n$ and 
\[
\B P(D^*_n \leq b) \ = \ \B E g((Z_1,\ldots,Z_n)|T_n=n) =
 \B E (R_b /n| T_n=n ) .
\]
With  $\epsilon,\rho,b,n$ and $q$ as above in (\ref{qrelation}), and
with the same $c(\rho)$ as above but with a different
$n(\epsilon,\rho)$, for all $n \geq n(\epsilon,\rho)$ we have
 $\B E U_b > \epsilon \ n
/\log n$ and $\B E R_b / n \ \in [q/\rho,q \rho]$.
Large deviation theory says that, with  $\lambda= \B E U_b$
as the mean of an unweighted sum of independent Poissons, the weighted
sum $Y=R_b$ satisfies  $\B P(Y/\B E Y \leq 1/\rho)
 \leq \exp(-\lambda c)$ and $\B P(Y/\B E Y \geq \rho)
 \leq \exp(-\lambda c)$. Putting
these together, we have for $n \geq n(\epsilon,\rho)$
\[
\B P(\frac{R_b}{n} \notin [q/\rho^2,q \rho^2] \ ) \leq
 2 \exp( -c(\rho) \epsilon \ n / \log n).
\]
Since the functional $g$ takes values in $[0,1]$, this proves, for $n
\geq n(\epsilon,\rho)$,
\eq  \label{biased with a rate}
\left| \B P(D^*_n \leq b) - q \right| \  \leq q(\rho^2 -1) +
 2 \exp( -c(\rho) \epsilon \ n / \log n)/\B P(T_n=n).
\en 

\subsubsection{The number of distinct block sizes}

 Odlyzko and Richmond (1985) prove that the number $J_n$ of distinct
block sizes in a random partition of the set $[n]$ is asymptotic to $e
\log n$ in expectation and in probability.  A stronger result can easily
be proved by overwhelming the conditioning.

Informally, our argument is that for $1 \leq i \leq (e-\epsilon)  \log
n$, the Poisson parameter $\lambda_i=x^i/i!$ is large, so that $\B P(Z_i=0)$ is
very small, in fact small enough to overwhelm the conditioning on
$\{T_n=n \}$, so that $\B P(C_i(n)=0)$ is also very small, and we
can conclude $\B P(C_i(n) = 0 $ for any $i \leq (e-\epsilon)  \log n))
\rightarrow 0$.  This accounts for at least $(e-\epsilon)  \log n$
distinct block sizes.  On the other side, $\sum_{i \geq (e+\epsilon) 
\log n} \B EZ_i$ is  small, hence
 for some $k=k(\epsilon)$, $\B P(Z_i >0$ for at least $k$ values
of $i \geq (e+\epsilon)  \log n)$ is very small, 
in fact small enough to overwhelm
the conditioning (using roughly $k=1/(2\epsilon)$.)
We conclude $\B P(C_i(n) > 0$ for at least $k$
 values $i \geq (e+\epsilon) \log n) \rightarrow 0$. Our result, that for any
$\epsilon > 0, \B P(C_i(n)=0$ for any $i \leq (e-\epsilon)  \log n$, or
$C_i(n)>0$ for at least $k$ values 
 $i \geq (e+\epsilon)  \log n) \rightarrow 0$,
implies but is not implied by the result that $J_n/\log n \rightarrow e$ in
probability. Furthermore, the bounds supplied by Theorem \ref{part
parts} below imply that $J_n/\log n \rightarrow e$ in $r$th mean for every $1
\leq r < \infty$. The result that $\B P(C_1(n)=0)\rightarrow 0$ was
proved in Sachkov (1974).

In a little more detail, observe that $\B P(Z_1=0)=\exp(-\lambda_1)=
e^{-x}=x/n \sim \log n /n$, which is smaller than the conditioning
probability, given by (\ref{Moser}), by a factor on the order of
$\sqrt{n}/(\log n)^{3/2}$.  The preceding argument is given in Sachkov
(1974).  The Poisson parameters increase rapidly, so
$\B P(Z_2=0)=\exp(-\lambda_2)=\exp(-x^2/2)=(x/n)^{x/2}$, which decays
faster than any power of $n$.

For a more careful analysis of the boundary where the Poisson parameter
$\lambda_i$ changes from large to small, write $i=(x+d)e$, where
$d=o(x)$. Recall $x \sim \log n$. 
 Using Stirling's formula, and writing $\approx$ for
 logarithmically asymptotic, we have $\lambda_i = x^i/i! \sim
(xe/i)^i/\sqrt{2 \pi i} = (x/(x+d))^i / \sqrt{2 \pi i} \approx $ $
\exp(-id/x -\log \sqrt{i}) \approx 
\exp(- ed -\frac{1}{2} \log \log n)$, so that the critical boundary for
$i$,
corresponding to
$d=-\frac{1}{2e} \log \log n$, is at $c(n) \equiv  xe - \frac{1}{2} \log
\log n$.  On the left side of this boundary 
the argument via overwhelming the
conditioning shows that $\B P(C_i(n)=0$ for any $i<ex-
(\frac{3}{2}+\epsilon) \log \log n) \rightarrow 0$. The
argument is very asymmetric between left and right: on the left, where
$\lambda_i$ is large, we use $\B P(Z_i=0)=\exp(-\lambda_i)$, gaining the
use of an exponential; while on the right, where $\lambda_i$ is small,
we use $\B P(Z_i>0)<\lambda_i$.  Thus in Theorem \ref{part parts}, the
left boundary $a$ is an extra $(1+\epsilon) \log \log n$ below $c(n)$,
while the right boundary $b$ is an extra $\epsilon \log n$ above
$c(n)$. 

The results of the above discussion are summarized by the following 
\begin{theorem}\label{part parts}
For partitions of a set of size $n$, for $\epsilon >0$, there are with
high probability blocks of every size $i \leq (e-\epsilon)  \log n$,
and not many blocks of size $i \geq (e+\epsilon)  \log n$.  More precisely,
for any $r<\infty$ there exists $k=k(\epsilon,r)<\infty$ so that, as $n
\rightarrow \infty$,
\[
\B P(C_1(n) = 0)  =  O((\log n)^{3/2}/\sqrt{n}),
\]
while for $a \equiv ex- (\frac{3}{2}+\epsilon) \log \log n)$
\[
\B P(C_i(n)=0 \mbox{ for any } 2 \leq i \leq a) \leq \frac{1}{\B P(T_n=n)} \sum_2^a
e^{-\lambda_i} = o(n^{-r}),
\]
and
\[
\B P(\sum_{i \geq b \equiv (e+\epsilon) \log n} C_i(n) \geq k)
=O \left(\frac{1}{\B P(T_n=n)}( \sum_{i\geq b} \lambda_i)^k \right)
=o(n^{-r}),
\]
where  $xe^x=n$, $\lambda_i=x^i/i!$,
and $\B P(T_n=n)$ satisfies (\ref{Moser}).
\end{theorem}
\proof
Most of the proof is contained in the informal discussion before the
theorem. For the second statement, it remains to check that $\sum_2^a \exp(-
\lambda_i)=o(n^{-r})$ for any $r$, which follows from an upper bound on
the first and last terms of the sum, which has at most $n$ terms, 
together with the observation that
the $\lambda_2 < \lambda_3 < \cdots < \lambda_{\lfloor x \rfloor} \geq
\cdots > \lambda_{\lfloor a \rfloor}$. For the third statement, we are
merely using the estimate, for $Y = \sum_{i \geq b} Z_i$, which is
Poisson with small parameter $\lambda$, that $\B P(Y \geq k)
=O(\lambda^k)$ as $\lambda \rightarrow 0$. Note that $\BE Y \approx \BE
Z_{\lceil b \rceil} \equiv \lambda_{\lceil b \rceil} \approx (xe/b)^b
\approx (1 + \epsilon/e)^{- b} < n^{-\epsilon}$.
\hfill \qed

The above argument by overwhelming the conditioning is crude but easy to
use because it gives away a factor of $\B P(T_n=n)$, when in fact the
event $\{T_n=n \}$ is approximately independent of the events involving
$\{Z_i>0\}$ for large $i$. An effective way to quantify and handle this
approximate independence is the total variation method outlined in
sections 3 and 4. Sachkov (1974) analyzed the size $L_n$ of the 
largest block of a random 
partition, and gave its approximate distribution.  Writing $L_n = h(\b
C(n))$ where $h(a_1,\ldots,a_n) = \max(i:a_i>0)$, Sachkov's result can
be paraphrased as $d_{TV}(L_n, h(\b Z_n)) \rightarrow 0$.  Note that the
number $J_n$ of distinct block sizes satisfies $J_n
\leq L_n$ always. 
Using $B=\{ i\leq n: i>ex - 2\log \log n \}$ for example, it should be
possible to prove that $d_{TV}(\b C_B,\b Z_B) \rightarrow 0$.  Then, by
comparison of $J_n = h(\b C(n))$ with $h(Z_1,\ldots,Z_n) = \sum 
\bone(Z_i>0)$, it would follow that, with centering constants $c(n)
\equiv ex -\frac{1}{2e} \log \log n$, the family of random variable
  $\{ J_n - c(n) \}$
is tight, and the family  $\{L_n-J_n\}$
is tight;  and for each family, 
along a subsequence $n(k)$ there is convergence in
distribution if and only if the  $c(n(k)) \mbox{ mod } 1$ converge.

 \section{Dependent process approximations}\label{sect12}

For the logarithmic class of structures discussed in Sections
\ref{logsect}, \ref{esfsect}, and \ref{sect5.2}, we have seen that the Ewens
sampling formula (ESF) plays a crucial role. In the counting process for
large components of
logarithmic combinatorial structures, there is substantial dependence;
an appropriate comparison object is the dependent process of large
components in the ESF. For example, in Arratia, Barbour
and Tavar\'e (1993) it is shown that the process of counts of factors of large
degree in a random polynomial over a finite field is
close in total variation to the process of counts of large cycles in a
random  permutation, corresponding to the ESF with parameter $\theta =
1$. In Arratia, Barbour and Tavar\'e (1994), Stein's method is used to
establish an analogous result for all the logarithmic class, and
somewhat more generally. The basic  technique involving Stein's method
specialized to the compound Poisson is described in Barbour, Holst
and Janson (1992, Chapter 10). 

Once such bounds are available, it is  a simple matter to establish
approximation results, with bounds, for other interesting functionals of
the large component counts of the combinatorial process. For example, 
the Poisson-Dirichlet and GEM limits for random polynomials are established with
metric bounds in Arratia, Barbour and Tavar\'e (1993). 
Poisson-Dirichlet limits for the logarithmic class are also discussed
by Hansen (1993).

\section{References}

\begin{verse}

\bigskip
APOSTOL, T.M. (1976) {\em An Introduction to Analytic Number Theory}, Springer
Verlag, New York.
 
ARRATIA, R. and TAVAR\'E, S. (1992a) The cycle structure of random
permutations. {\em Ann. Prob.} 20, 1567-1591.

ARRATIA, R. and TAVAR\'E, S. (1992b) Limit theorems for combinatorial
structures via discrete process approximations. {\em Rand. Struct. Alg.}
 3, 321-345.

ARRATIA, R., BARBOUR, A.D. and TAVAR\'E, S. (1992) Poisson process
approximations for the Ewens Sampling Formula. {\em Ann. Appl. Prob.}
2, 519-535.

ARRATIA, R., BARBOUR, A.D.  and TAVAR\'E, S. (1993)
On random polynomials over finite fields. {\em Math. Proc. Camb. Phil.
Soc.} 114, 347-368.

ARRATIA, R., BARBOUR, A.D.  and TAVAR\'E, S. (1994) 
Logarithmic combinatorial structures. In preparation.

ARRATIA, R., STARK, D. and TAVAR\'E, S. (1994) Total variation asymptotics 
for Poisson process approximations of logarithmic combinatorial 
assemblies. {\em Ann. Probab.}, to appear.

BARBOUR, A.D. (1992) Refined approximations for the
Ewens sampling formula.  {\em Rand. Struct. Alg.} 3, 267-276.

BARBOUR, A.D., HOLST, L., and JANSON, S. (1992) {\em Poisson
Approximation}. Oxford University Press, Oxford.

BOLLOB\'AS, B. (1985) {\em Random Graphs.} Academic Press, New York.

de BRUIJN, N.G. (1981) {\em Asymptotic Methods in Analysis}, Dover
Press. (Republication of 1958 edition, North-Holland Publishing Co.)

DeLAURENTIS, J.M. and PITTEL, B.G. (1983) Counting subsets of the random
partition and the `Brownian bridge' process.  {\em Stoch. Proc. Appl.}
 15, 155-167.

DIACONIS, P. and FREEDMAN, D. (1980) Finite exchangeable sequences.
{\em Ann. Prob.} 8, 745-764.

DIACONIS, P. and PITMAN, J.W. (1986) Unpublished lecture notes,
Statistics Department, University of California, Berkeley.

DIACONIS, P., McGRATH, M., and PITMAN, J.W. (1994) 
Riffle shuffles, cycles and descents.
{\em Combinatorica}, in press.

ELLIS, R.S. (1985) {\em Entropy, Large Deviations, and 
 Statistical Mechanics.} Springer, Berlin.
 
ETHIER, S.N. and KURTZ, T.G. (1986).  {\em Markov
Processes: Characterization and Convergence}, Wiley, New
York.

EWENS, W.J. (1972).  The sampling theory of selectively
neutral alleles.  {\em Theor. Pop. Biol.}  3, 87-112.

FELLER, W. (1945) The fundamental limit theorems in probability. {\em
Bull. Amer. Math. Soc.} 51, 800-832.

FLAJOLET, P. and ODLYZKO, A.M. (1990a) Singularity analysis of generating
functions. {\em SIAM J. Disc. Math.} 3, 216-240.

FLAJOLET, P. and ODLYZKO, A.M. (1990b) Random mapping statistics.
In {\em Proc. Eurocrypt '89}, J.-J. Quisquater, editor, pp. 329-354.
	Lecture Notes in C.S. 434, Springer-Verlag.

FLAJOLET, P. and SORIA, M. (1990) Gaussian limiting distributions for
the number of components in combinatorial structures. 
{\em J. Comb. Th. A} 53, 165-182.

FOATA, D. (1974) La s\'erie g\'en\'eratrice exponentielle dans les
probl\'emes d'\'enum\'erations. Press Univ. Montreal.

FRISTEDT, B. (1992) The structure of random partitions of large sets.
Preprint.

FRISTEDT, B. (1993) The structure of random partitions of large integers.
{\em Trans. Amer. Math. Soc.} 337, 703-735.

GOH, W.M.Y. and SCHMUTZ, E. (1993) The number of distinct
parts in a random integer partition. {\em J. Comb. Theory A}, in press.

GONCHAROV, V.L. (1944) Some facts from combinatorics. {\em Izvestia Akad.
Nauk. SSSR, Ser. Mat.} 8, 3-48. See also: On the field of combinatory 
analysis. {\em Translations Amer. Math. Soc.} 19, 1-46.

GRIFFITHS, R.C. (1988) On the distribution of points in a Poisson--
Dirichlet process. {\em J. Appl. Prob.} 25, 336-345.

HANSEN, J.C. (1989) A functional central limit theorem for random 
mappings. {\em Ann. Prob.} 17, 317-332.

HANSEN, J.C. (1990) A functional central limit theorem for the Ewens
Sampling Formula. {\em J. Appl. Prob.} 27, 28-43.

HANSEN, J.C. (1993) Order statistics for decomposable combinatorial
structures. {\em Rand. Struct. Alg.}, in press.

HANSEN, J.C. and SCHMUTZ, E. (1994) How random is the characteristic
polynomial of a random matrix? {\em Math. Proc. Camb. Phil.
Soc.} 114, 507-515.

HARPER, L.H. (1967) Stirling behavior is asymptotically normal. {\em
Ann. Math. Stat.} 38, 410-414.

HARRIS, B. (1960) Probability distributions related to random mappings.
{\em Ann. Math. Stat.} 31, 1045-1062.


HOLST, L. (1979a) A unified approach to limit theorems for urn models.
{\em J. Appl. Prob.} 16, 154-162.

HOLST, L. (1979b) Two conditional limit theorems with applications. {\em
Ann. Stat.}  7, 551-557.

HOLST, L. (1981) Some conditional limit theorems in exponential families.
{\em Ann. Prob.} 9, 818-830.

IGNATOV, T. (1982) On a constant arising in the asymptotic theory of symmetric
groups, and on Poisson--Dirichlet measures. {\em Theory Prob. Applns.}
27, 136-147.

JOYAL, A. (1981) Une th\'eorie combinatoire des s\'eries formelles.
{\em Adv. Math.} 42, 1-82.

KOLCHIN, V.F. (1976) A problem of the allocation of particles in cells
and random mappings. {Theor. Probab. Applns.} 21, 48-63.

KOLCHIN, V.F. (1986) {\em Random Mappings}, Optimization Software, Inc.,
New York.

KOLCHIN, V.F., SEVAST'YANOV, B.A., and CHISTYAKOV, V.P. (1978) {\em
Random Allocations.} Wiley, New York.

LEVIN, B. (1981) A representation for multinomial cumulative
distribution functions. {\em Ann. Stat.} 9, 1123-1126.

LIDL, R. and NIEDERREITER, H. (1986) {\em Introduction to Finite Fields
and their Applications}, Cambridge University Press.

MEIR, A. and MOON, J.W. (1984) On random mapping patterns. {\em
Combinatorica} 4, 61-70.

METROPOLIS, N. and ROTA, G.-C. (1983) Witt vectors and the algebra of
necklaces. {\em Adv. Math.} 50, 95-125.

METROPOLIS, N. and ROTA, G.-C. (1984) The cyclotomic identity. Contemporary 
Mathematics, Volume 34, 19 - 27. American Mathematical Society, Providence,
R.I.

MOSER, L. and WYMAN, M. (1955) An asymptotic formula for the Bell
numbers. {\em Trans. Roy. Soc. Canad.} 49, 49-53.

MUTAFCIEV, L.R. (1988) Limit theorems for random mapping patterns. {\em
Combinatorica} 8, 345-356.

ODLYZKO, A.M. and RICHMOND, L.B. (1985)  On the number of distinct block
sizes in partitions of a set.  {\em J.  Comb. Theory A} 38, 170-181.

OTTER, R. (1948) The number of trees. {\em Ann. Math.} 49, 583-599.

POURAHMADI, M. (1984) Taylor expansion of $\exp(\sum_{k=0}^\infty a_k
z^k)$ and some applications. {\em Amer. Math. Monthly} 91, 303-307.


SACHKOV, V.N. (1974) Random partitions of sets. {\em Th. Prob. Appl.} 19,
184-190.

SHEPP, L.A. and LLOYD, S.P. (1966) Ordered cycle lengths in a
random permutation. {\em Trans. Amer. Math. Soc.} 121, 340-357.

STAM, A.J. (1978) Distance between sampling with  and without
replacement. {\em Statistica Neerlandica} 32, 81-91.

STARK, D. (1994a) Unpublished Ph.D. Thesis, Department of Mathematics, 
University of Southern California.

STARK, D. (1994b) Total variation asymptotics for independent process
approximations of logarithmic multisets and selections. Preprint.

STEPANOV, V.E. (1969) Limit distributions for certain characteristics of
random mappings. {\em Theory Probab. Appl.} 14, 612-626.

VERSHIK, A.M. and SHMIDT, A.A. (1977)
Limit measures arising in the theory of groups I.
{\em Theor. Prob. Applns.} 22, 79-85.


WATTERSON, G.A. (1974a) The sampling theory of selectively neutral
alleles. {\em Adv. Appl. Prob.} 6, 463-488.

WATTERSON, G.A. (1974b) Models for the logarithmic species abundance
distributions. {\em Theor. Popn. Biol.} 6, 217-250.

WATTERSON, G.A. (1976) The stationary distribution of the infinitely
many alleles diffusion model. {\em J. Appl. Prob.} 13, 639-651.

WILF, H.S. (1990) {\em Generatingfunctionology}. Academic Press, San
 Diego, CA.

\end{verse}

 \end{document}